\documentclass[10pt]{amsart}
\usepackage{amssymb}
\usepackage{latexsym}
\usepackage{epsfig}
\usepackage{graphicx}
\usepackage{psfrag}
\usepackage{ amssymb, amsmath}
\usepackage{amsfonts}
\usepackage{amsthm}
\usepackage{enumerate}
\usepackage{color}
\usepackage{ifpdf}
\newtheorem{theorem}{Theorem}[section]
\newtheorem{corollary}[theorem]{Corollary}
\newtheorem{lemma}[theorem]{Lemma}
\newtheorem{proposition}[theorem]{Proposition}

\newtheorem{definition}[theorem]{Definition}

\def\be{\begin{equation}}
\def\ee{\end{equation}}
\def\bes{\begin{equation*}}
\def\ees{\end{equation*}}
\def\bea{\begin{eqnarray}}
\def\eea{\end{eqnarray}}
\def\eeas{\end{eqnarray*}}
\def\beas{\begin{eqnarray*}}
\def\supp{\hbox{supp}\,}

\def\coi{C_0^\infty}
\def\D{\Delta}
\def\dim{\textrm{ dim }}
\def\e{\epsilon}
\def\hra{\hookrightarrow}
\def\iohoh{I^{-\frac12,\frac12}}

\def\Ker{\textrm { Ker }}
\def\la{\langle}
\def\lra{\longrightarrow}
\def\o{\omega}
\def\p{\partial}
\def\pl{\pi_L}
\def\pr{\pi_R}
\def\ra{\rangle}

\def\R{\mathbb R}

\def\tC{{C}}

\def\tx{T^*X\setminus 0}
\def\ty{T^*Y\setminus 0}

\def\ajk{A_{jk}}

\def\dim{\textrm{ dim }}

\def\Ker{\textrm { Ker }}
\def\lra{\longrightarrow}

\def\R{\mathbb R}

\def\tho{\theta_1}
\def\tht{\theta_2}

\def\at{\tilde\alpha}
\def\bt{\tilde\beta}
\def\rt{\tilde\rho}
\def\st{\tilde\sigma}
\def\xt{\tilde\xi}

\def\agtxt{ }

\def\tP{\tilde\Phi}

\numberwithin{equation}{section}

\begin{document}

\title[An FIO Calculus: Sobolev Estimates]{An FIO calculus for marine seismic imaging, \\
 II:  Sobolev estimates}

\author{Raluca Felea,  Allan Greenleaf and Malabika Pramanik}
\thanks{The second author was partially supported by NSF grants
DMS-0138167 and DMS-0551894. The third author was partially supported by an NSERC Discovery grant.}

\date{Revised, Feb. 7, 2010}

\maketitle

\begin{abstract}
We establish sharp $L^2$-Sobolev estimates for classes of pseudodifferential operators with singular symbols \cite{meuh,guuh} whose non-pseudodifferential (Fourier
integral operator) parts exhibit two-sided fold singularities. The operators
considered include both singular integral operators along curves in $\R^2$ with simple inflection points and normal operators arising in linearized seismic
imaging in the presence of fold caustics \cite{no, fel, FG}. 
\end{abstract}

\section{Introduction}\label{intro}
{\allowdisplaybreaks{
This paper is concerned with $L^2$-Sobolev estimates for some classes of operators sharing a common microlocal geometry that displays two different kinds of singularities. One of our results is that in certain cases these singularities interact, producing worse estimates than one might otherwise expect. We begin by describing the main features of this geometry and the motivations for studying these operators.  

\subsection{The cubic model} \label{subsec-cubic} In $\R^2$,  consider the Hilbert transform along a model curve $\gamma(t):=(t,t^3)$ with a simple inflection point,
\begin{equation}\label{eqn-hil}
\mathcal Hf(x)= p.v. \int_\R f(x_1-t,x_2-t^3) \frac{dt}t.
\end{equation}
It is  well known  that  $\mathcal H$ (and similar but much more general
operators) are bounded
on $L^2(\R^2)$ (see, e.g., \cite{SW1,SW2,CNSW}.) The Schwartz kernel of $\mathcal H$ is given by
\begin{equation}
K_{\mathcal H}(x,y)=\delta(x_2-y_2-(x_1-y_1)^3)\times p.v.\left(\frac1{x_1-y_1} \right),
\end{equation}
which incorporates both pseudodifferential and Radon transform type
singularities. In fact, the wave front set of  $K_{\mathcal H}$ satisfies the inclusion $WF(K_{\mathcal H})\subset \D'\cup \tC_0'$, where \begin{itemize} \item $\D = \{(x, \xi; x, \xi) : x \in \mathbb R^2, \xi \in \mathbb R^2 \setminus \{0\} \}$ is the diagonal of $T^*\R^2\times T^*\R^2$, \item $\tC_0'$ is the conormal bundle of $\{x_2-y_2-(x_1-y_1)^3=0\}$, and 
\item $\, '\, $ denotes the
usual twist map in microlocal analysis, sending \linebreak$(x,\xi;y,\eta)\mapsto
(x,\xi; y,-\eta)$. \end{itemize}The same geometry is present for fractional integral operators along $\gamma$, 
\begin{equation}\label{def-frac}
\mathcal J_{\alpha} f(x)=\int f(x_1-t,x_2-t^3)\frac{dt}{|t|^\alpha},\quad 0<\alpha<1.
\end{equation}
Similar operators, combining fractional integral singularities along the curve with fractional derivatives transverse to it,  have been studied for various curves in $\R^2$ or in nilpotent groups \cite{RS,Graf,CS}. 
Operators with kernels that display two different types of singularities arise in a variety of settings and go under various names: In the terminology of \cite{PS1}, $\mathcal H$ is a \emph{singular Radon transform}, while in the language of \cite{meuh,guuh},
$\mathcal H$ and $\mathcal J_{\alpha}$ are \emph{pseudodifferential operators with
singular symbols}. The latter concept was originally introduced to describe parametrices for  variable coefficient wave and other real
principal type operators, but has since found numerous applications in integral geometry \cite{gruh2,gruh3} and  inverse scattering \cite{GRUHCMP,KLU}. 

{\agtxt In the present paper we obtain sharp $L^2$-Sobolev estimates for a class of operators containing $\mathcal H$ and $\mathcal J_{\alpha}$, and with similar  microlocal geometry. 
One feature is that in some cases the estimates  are worse than one would expect by considering the orders of the pseudodifferential and Fourier integral operator parts of the operators separately.
For example,  $\mathcal J_{\alpha}:H^s_{\text{comp}}\lra H^{s-r}_{\text{loc}}$ for $r=\frac{\alpha-1}3$, rather than for  $r=\max(\alpha-1,-\frac13)$. Such behavior can already be seen, for general families of curves (for $\alpha$ close to 1) in Greenblatt \cite{mG}.

However,  the Schwartz kernels of the operators we are most interested in are not associated with a nested pair of submanifolds  \cite{gruh2} in the same way that those of fractional integrals along submanifolds are. We thus  need to formulate this class of operators using a microlocal approach; from this point of view,} $\mathcal H\in\iohoh(\Delta,\tC_0)$ and $\mathcal J_\alpha\in I^{-\frac12,\alpha-\frac12}(\Delta,\tC_0)$, where, for $p,l\in\R$, the space
$I^{p,l}(\D,\tC_0)$ denotes the class of operators on
$\mathcal E'(\R^2)$ whose Schwartz kernels belong to the class of
paired Lagrangian distributions associated with the cleanly
intersecting Lagrangian manifolds $\D'$ and ${\tC_0}'$ in $T^*(\R^2\times\R^2)$. We refer the reader to \cite{meuh,guuh,men,gruh2} for  the theory of the $I^{p, l}$ classes and to \S \ref{sec-back} for the basic definitions that we need here. The key facts relevant for this paper are that 
\begin{enumerate}[(i)] 
 \item {\agtxt  The canonical relations  $\D,\, {\tC_0}$ intersect cleanly in codimension one \nolinebreak\cite{gruh2}.}
  \item  $\tC_0$ is a \emph{folding} canonical relation in the sense of Melrose and Taylor\cite{meta} or a \emph{two-sided fold} in the sense of \cite{grse,grse2} (see \S\S\ref{prelim} for details). 
  \item The fold surface of $\tC_0$ equals its intersection with the diagonal,  $\D\cap\tC_0$. 
  \end{enumerate}
  {\agtxt We consider the pair $(\Delta, \tC_0)$ mainly as a stepping-stone to the  seismic geometries described below, and  do not study non-translation invariant versions of $\mathcal H$ and $\mathcal J_\alpha$  along the lines of \cite{RS,Graf,se93,CS}.}

\subsection{Single-source Seismic Imaging} \label{subsec-single} A very similar geometry occurs in a completely unrelated problem, involving forward scattering maps  in   linearized seismic imaging \cite{no,fel,FG}.  Let $Y=\R^3_+:=\{y_3>0\}$ model the subsurface of the earth, and $\partial Y=\{y_3=0\}$ its surface or the surface of the ocean. Let $\Sigma\subset \partial Y\times \partial Y$ denote a \emph{source-receiver manifold} of  pairs $(s,r)$, and $T >0$ the total time-length of the seismic experiment. The \emph{data space} for a typical such imaging experiment is given by the collection $X=\Sigma\times (0,T)$, which corresponds to surface measurements being made at locations $r$ and times $t$, $0 < t < T$,
resulting from seismic events with (idealized) delta-function impulses at  sources $s $ and time $t_0=0$. The dimension of the data space is of course one more than that of the source-receiver manifold. It is known that under fairly general conditions, the corresponding (formally) linearized scattering operator $F$, mapping perturbations in the sound speed $c(y)$ in $Y$ (about a known smooth background $c_0(y)$) to perturbations in the pressure field measured on $X$, is a Fourier integral operator (hereafter referred to as FIO)  associated with a canonical relation $C\subset (T^*X\setminus 0)\times (T^*Y\setminus 0)$ \cite{ra,tKSV,nosy}. More specifically, for some order $m \in \mathbb R$ dictated by the dimensionality of the problem at hand, $F$ belongs to the  class  of FIOs, $I^m(X, Y; C)$, consisting of operators mapping $\mathcal E'(Y) \longrightarrow \mathcal D'(X)$ whose Schwartz kernels are Fourier integral distributions associated with $C'\subset T^*(X\times Y)\setminus 0$. 

The $L^2$-Sobolev regularity  of the operator $F$, which in turn is tied to the question of microlocal invertibility, turns out to be heavily dependent on the data collection geometry and certain non-degeneracy hypotheses regarding the background sound speed $c_0(y)$, assumed smooth  and known. Historically, the first situation to be mathematically analyzed was the {\em{single source geometry without caustics}}, for which (a) the seismic data is obtained from a single source, $s_0\in\partial Y$ and  the receivers range over an open subset of $\partial Y$ (i.e., $\dim X = \dim Y = 3$); and (b)  the  rays flowing out from the source $s_0$ have no caustics. The latter assumption means the following. Let us consider the smooth conic Lagrangian $\Lambda_{s_0}\subset T^*\R^3\setminus 0$ obtained from the flowout of the points $(\xi_0, \tau_0)$ in the characteristic variety of the wave equation above the point $(x,t)=(s_0, 0)$ using the Hamiltonian vector field of the principal symbol of the background wave operator, $\square_{c_0}$. There will be no caustics if the spatial projection from this Lagrangian manifold, $\pi_Y:\Lambda_{s_0}\longrightarrow Y$ , has full rank everywhere (except of course $(s_0,0)$). Beylkin \cite{be} showed that under these assumptions $C$ is a local canonical graph, so that $F^*F$ is an elliptic pseudo-differential operator (hereafter referred to as $\Psi$DO) on $Y$. Application of a left parametrix of $F^*F$ then results in \emph{high-frequency linearized seismic inversion}. This conclusion also holds for other data geometries, under the \emph{traveltime injectivity condition} \cite{tKSV,nosy,stolk}, which ensures that the canonical relation $C$ satisfies the Bolker condition \cite{Gu}, i.e., $\pi_L:C\longrightarrow T^*X\setminus 0$ is an injective immersion; this was further weakened in some situations \cite{stolk}.

However, caustics are unavoidable in physically realistic
velocity models  \cite{nosy-fold} and, since they typically lead to artifacts in
images, understanding the
structure of $F$ and
$F^*F$ when the background sound speed exhibits caustics is a
fundamental problem in
exploration seismology.  Nolan
\cite{no} showed that, for the single source geometry with caustics of at most fold type, $C$ is a folding canonical relation.
It was then shown  that, for \emph{any} folding relation $C_{fold}$ and any $F\in I^{m}(C_{fold})$,
\be\label{ss orders}
F^*F\in I^{2m,0}(\D,C_1),
\ee
 with $C_1$ having the  properties  (i,ii,iii) as in \S\S\ref{subsec-cubic} above  (\cite{no},  Felea \cite{fel}). In particular, for the seismic problem, $F^*F\in I^{2,0}(\D,C_1)$ .

\subsection{Marine seismic imaging}\label{subsec-marine}
Another seismic data set of interest comes from  \emph{marine} (or {offshore})  imaging. A mathematical idealization of the experimental setup is as follows. A survey vessel trails behind it a cable containing both an acoustic source  and a line of recording instruments. The point source consists of an airgun which sends acoustic waves through the ocean to the subsurface. Reflections occur when the sound waves encounter singularities  in the soundspeed in the subsurface, such as discontinuities at interfaces of sedimentary layers, and the reflected rays are then received by a linear array of hydrophones towed behind the vessel. The vessel  then makes repeated passes along parallel paths (say, parallel to the $x_1$ axis) contained in an open  $U\subset\partial Y$. In other words, the source-receiver pairs form an open subset of
\[ \Sigma = \{(r, s) \in U \times U : r= (r_1, r_2, 0), s = (s_1, s_2, 0), r_2 = s_2 \}, \]
which is a codimension one in $\partial Y\times \partial Y$. Thus the
data set is overdetermined, with $\dim X=4 > \dim Y=3$. The forward operator $F\in I^{1-\frac14}(X,Y;C)$, and in \cite{FG} the first two authors  identified the structure of $C$, under a natural extension of the fold caustic assumption to this  context. General canonical relations having this structure, along with some additional  nondegeneracy conditions, were  called \emph{folded cross caps}, and a composition result was proven: if  $C$ is a folded cross cap and $A \in I^{m-\frac{1}{4}}(X,Y;C)$, then 
\be\label{marine orders}
A^*A \in I^{2m- \frac{1}{2},\frac{1}{2}}(\Delta, C_2),
\ee
with  $C_2 \subset \left(T^*Y \setminus 0\right)
\times \left(T^*Y \setminus 0\right)$  again satisfying (i,ii,iii) from \S\S\ref{subsec-cubic}.  In particular, for the marine seismic imaging problem, $F^*F\in I^{\frac32,\frac12}(\Delta, C_2)$.

Now we can compare the normal operators for the single source and  the marine geometries (both in the presence of fold caustics). It follows from \cite{meuh,guuh} that microlocally away from $\D\cap C_2$ and for any $p, l \in \mathbb R$, 
\begin{equation}\label{away}
I^{p,l}(\D,C_2)\hra I^{p+l}(\D\setminus C_2) +
I^{p}(C_2\setminus\D).
\end{equation}
From (\ref{marine orders}) and (\ref{away}), one observes that since $l = \frac{1}{2}$ for the marine geometry, the order of the non-$\Psi$DO part,
i.e., the FIO part, of the normal operator $F^*F$ is $\frac{1}{2}$ lower than the $\Psi$DO part. On the other hand, from  (\ref{ss orders}) and the analogue of  (\ref{away}) one sees that for the single source geometry the two orders are the
same. In other words, the artifact resulting from the FIO part of the normal operator for the single source data set is at least as strong as the $\Psi$DO part and hence is nonremovable, but the fact that the artifact for the marine data set is $\frac12$ lower order (away from $\Delta\cap C_2$) leads one to hope that in this situation  the artifact might be removable.

\subsection{Objectives and scope} However, the explanation in the previous paragraph is informal in the sense that it ignores the singular behavior of the operators at the intersection 
$\D\cap C_j,\, j=0,1,2$, and the intuition
above needs to be justified. As an initial step, one should understand the mapping properties of operators in $I^{p,l}(\D,C_j)$.  The primary goal of this paper is to obtain sharp $L^2$ Sobolev estimates for operators in $I^{p,l}(\D,C_j)$, $j=0,1,2$.

Although there exists a microlocal normal form for a folding canonical relation, $C$, obtained by applying suitable canonical transformations on
the left and right \cite{meta}, when applied to $\D$ these same transformations will turn $\D$ into some unknown canonical graph. Thus, there is no such normal form for  pairs such as
$(\D,C)$ satisfying (i,ii,iii) in \S\S\ref{subsec-cubic}. For this reason, our results will be limited to the pairs $(\D, C_j),\, j=0,1,2$, described above. In other words, our results apply to  general
operators with the same cubic  geometry underlying $\mathcal H, \mathcal J_\alpha$; the geometry arising from $F^*F$ for any FIO with a folding canonical relation, including the  single source forward operator; and  the geometry for the normal operator in the  marine seismic imaging problem. In these cases, one has
good control over the phase functions  in the oscillatory
representations of the operators. 

\subsection{Main result} We show 
that, for the canonical relations that we consider, the following Sobolev estimates hold:
\begin{theorem}\label{thm-main}
Let $A\in I^{p,l}(\D,C_0)$ or $A\in I^{p,l}(\D,C_1)$ or $A\in I^{p,l}(\D,C_2)$, where the canonical relations $C_j$ are as above. Then $A:H^s_{\text{comp}}\lra H^{s-r}_{\text{loc}}$ for
\begin{equation}\label{r-in-main}
r=\begin{cases} p+\frac16, & l<-\frac12,\\
p+\frac16+\epsilon, & l=-\frac12,\, \forall \e>0,\\
p+\frac{l+1}3, & -\frac12<l<\frac12,\\
p+l, & l\ge \frac12.
\end{cases}
\end{equation}
\end{theorem}

\noindent{\bf Remarks:} \begin{enumerate}[1.]
\item Here $H^s$,  $H^{s}_{\text{comp}}$ and $H^{s}_{\text{loc}}$  denote the standard $L^2$-based Sobolev space and its compactly supported and local variants.
Thus our analysis of the Schwartz kernels can be restricted to compact sets in the spatial variables.   
\medskip 

\item Away from $\D\cap C_i$, which is also the fold surface of $C_i$,
the inclusion (\ref{away})
holds, and $I^{p+l}(\Delta \setminus
C_i)$ and $I^p(C_i \setminus \Delta)$ satisfy standard Sobolev space
estimates associated with local canonical graphs. These estimates are always
more regular than the ones mentioned in the statement of Thm.\
\ref{thm-main}. Thus it suffices to only consider operators $A$
supported microlocally close to $\Delta \cap C_i$. Keeping in mind
both the $\Psi$DO and FIO natures of $I^{p,l}(\Delta, C_i)$ and the
loss of $\frac16$ derivative for FIOs associated with folding canonical relations
\cite{meta},  in general one certainly needs to take
$r\ge\max(p+l, p+\frac16)$. Thm.\ \ref{thm-main} can therefore be interpreted as saying that, when the strengths of the $\Psi$DO and FIO parts of $A$ are
sufficiently close, specifically $-\frac12\le l<\frac12$, there is a further loss due to their interaction. We will
see in \S\ref{sec-sharp} that this loss can in fact occur and Thm.\ \ref{thm-main} is in general sharp. 
For $l<\frac12$ close to $\frac12$, this type of behavior is already present in  the estimates of Greenblatt \cite{mG} for fractional integrals along families of curves.
\medskip

\item For comparison, it is  natural to investigate the analogue of Thm. \ref{thm-main} when $C$ is a local canonical graph intersecting $\Delta$ cleanly in codimension one. We address this in Thm. \ref{thm-cg}. Singular Radon transforms over hypersurfaces
satisfying the rotational curvature condition
are of this general type but with higher codimension intersections.  A particular case was first considered by Geller and Stein \cite{GeSt} and the general class is due to Phong and Stein \cite{PS1}; see also Cuccagna \cite{CuThesis}, which allows for degenerate canonical relations.  Different proofs of $L^2$-boundedness   were given in \cite{gruh2, PS2}. Our proof of Thm. \ref{thm-cg} is essentially an adaptation of arguments in \cite{gruh2} using the methodology of parabolic cutoffs, cf. \cite{mel}, whereas the more degenerate geometry of Thm.\ \ref{thm-main} requires a combination of   a parabolic-type decomposition with a two-index dyadic decomposition inspired by the almost orthogonal decompositions of degenerate oscillatory integral operators  due to  Phong and Stein \cite{PS3,PS4}.      
\end{enumerate}

\subsection{Layout of the paper}\label{layout}
 In \S \ref{sec-back}, we recall the basic theory of Fourier integral
operators associated with a single smooth canonical relation, as well
as the paired-Lagrangian operators associated with two cleanly
intersecting canonical relations. The oscillatory representation of
the two classes of operators under consideration, in particular 
the normal operators for the linearized marine seismic imaging problem
in the presence of fold caustics \cite{FG} is also reviewed. In \S \ref{sec-sharp}, we prove that the regularity exponents in Thm.\ \ref{thm-main} are sharp, in particular showing that the loss of derivatives for $I^{p,l}(\D, C_0)$ in the range $-\frac12\le l\le \frac12$ in Thm.\ref{thm-main} is unavoidable. As a warm-up to the proof of Thm.\ref{thm-main}, in \S \ref{sec-graph} we
give an analogous result in the case when the canonical relation is a
canonical graph; here, a standard parabolic decomposition suffices to
establish the Sobolev boundedness. In \S \ref{sec-gendecomp} and \S
\ref{sec-norm} we introduce, for the classes $I^{p,l}(\Delta, C_0)$
and $I^{p,l}(\Delta, C_2)$ respectively, the phase space decomposition of mixed
type which is the main technical feature of the paper. The proof of
Thm.\ \ref{thm-main} is given in \S \ref{sec-cubic}. Finally, in \S
\ref{hilbert-example}, we examine the Hilbert transform $\mathcal H$ along the cubic in more detail. We show that, although elliptic in a naive sense, $\mathcal H$ is not microlocally invertible using operators bounded on $L^2$, indicating the difficulty of constructing parametrices for these kinds of operators.

We thank the referee for suggestions concerning  exposition and references.

\section{Background} \label{sec-back}
 In this section we develop the notation and terminology needed in the sequel, and in particular provide an explicit representation for the operators of interest. We refer the reader to  \cite{ho, du,ho-book} for the theory of classical FIOs, 
\cite{meuh,guuh,men,gruh2} for material on the $I^{p,l}$ classes and parabolic cutoffs, and \cite[\S3]{FG} and the references there for further details on the relevance of these techniques for linearized seismic imaging. However, for convenience we briefly recall the salient facts in the form  needed here.  

\subsection{Preliminaries} \label{prelim} Let $X$ and $Y$ be  smooth manifolds and  $\tx,\ty$ their cotangent bundles (with the zero-sections deleted), equipped with the canonical symplectic forms $\,\o_{T^*Y},\,\o_{T^*Y}$.
If $\Lambda\subset
(\tx)\times(\ty)$ is a conic Lagrangian submanifold with respect to
$\o_{T^*(X\times Y)}=\pl^*\o_{T^*X}+\pr^*\o_{T^*Y}$, then
$\Lambda':=\{(x,\xi;y,\eta): (x,\xi;y, -\eta)\in\Lambda\}$ is
called a {\em{canonical relation}}. For example, the conormal bundle of any smooth submanifold $Z\subset Y\times Y$ is a canonical relation.  In particular,
if $X=Y$, then $\D=N^*\D_Y'\subset (\ty)\times(\ty)$ is the diagonal relation.

If $\Lambda$ is a Lagrangian,  $m \in\R$ and $0 \leq \delta\leq 1-\rho \leq \rho \leq 1$, then $I^m_{\rho, \delta}(X \times Y, \Lambda)$ denotes the space of Fourier integral distributions of order $m$ and type $(\rho, \delta)$ whose wavefront set is contained in $\Lambda$ (see \cite[p. 97]{du}), while if $C$ is a canonical relation, the associated class $I^m_{\rho,\delta}(X,Y;C)$ of FIOs of order $m$ and type $(\rho, \delta)$ 
mapping $\mathcal E'(Y) \lra \mathcal D'(X)$ consists of those operators $A$  whose Schwartz
kernels are Fourier integral distributions $K_A\in I^m_{\rho, \delta}(X\times Y;C')$ \cite{ho,du}. This is  abbreviated to $I^m_{\rho, \delta}(C)$ if $X$ and $Y$ are understood,
and if $\rho = 1$ and $\delta = 0$, the subscripts $\rho$ and $\delta$ are omitted. If $\pl,\pr:C\lra\ty$ are the natural projections on the left and right, then a canonical
relation $C$ is a \emph{local canonical graph} if and only if one of the
projections (and hence the other) is a local diffeomorphism; of course, this can only happen if $ \dim X=\dim Y$. In this case, $A\in I^m(C)\implies A:H^s_{comp}(Y)\lra H^{s-m}_{loc}(Y)$. In particular, if $C = \Delta$, then $I^m_{\rho,\delta}(\Delta)$ is the class of pseudodifferential operators of order $m$ and type $(\rho,\delta)$, $\Psi^m_{\rho,\delta}(Y)$. 

If $C$ fails to be a canonical graph, then the critical points of the
two projections are the same, and $\dim\Ker(D\pl)=\dim\Ker(D\pr)$ at
all points. One says that $C$ is a
\emph{folding canonical relation} \cite{meta}, also called a
\emph{two-sided fold} \cite{grse}, if both projections have only Whitney fold singularities; in this equidimensional setting, that is the same as $S_{1,0}$ singularity in the Thom-Boardman notation \cite{gogu}. A well known
result of Melrose and Taylor \cite{meta} states that, for such a $C$, there is a loss of $\frac16$ derivative: 
\be\label{est meta}
 A\in I^m(C)\implies A:H^s_{comp}\lra H^{s-m-\frac16},\, \forall s\in\R. 
 \ee  
{\bf{Remark.}} Note that, combined with the composition result of \cite{no,fel}, this gives an example exhibiting the loss in Thm. \ref{thm-main} for $l=0$. In fact, let $C$ be the folding canonical relation in the single source fold-caustic seismic geometry 
\cite{no}, and let $A\in I^m(C)$ be properly supported. By (\ref{ss orders}) the operator $A^*A\in I^{2m,0}(\D,C_1)$, with $(\D, C_1)$ as in Thm. A. Since the loss of $\frac16$ derivative for  $A $ is in general sharp, so is $A^*A:H^s\lra H^{s-2m-\frac13}$, already showing the necessity  of the loss of $\frac{l+1}3$ derivative for $l=0$.

\subsection{Paired Lagrangian spaces} 
For an exposition on clean intersection theory and Fourier integrals, we refer the reader to \cite{gogu, du}. Classes of distributions associated with two cleanly intersecting Lagrangian manifolds were introduced by Melrose and Uhlmann \cite{meuh} and Guillemin and Uhlmann \cite{guuh} (see also \cite{gruh2} for further relevant discussion). For our purposes, the definition of this class can be given in terms of multiphase functions \cite{men} and symbol-valued symbols \cite{guuh,meuh}, which we now describe.  

\begin{definition}
Given a manifold $X$ and a positive integer $M \geq 1$, let $\Gamma$ be a cone in $X \times \mathbb R^M \setminus \{0\}$. 
\begin{enumerate}[(i)]
\item A function $\varphi = \varphi(x; \theta) \in C^{\infty}(\Gamma)$ is a {\em{phase function}} if it is homogeneous of degree 1 in $\theta$ and has no critical points as a function of $(x, \theta)$.   
\item A phase function $\varphi$ is said to be {\em{non-degenerate in $\Gamma$}} if 
\begin{align*} 
&d_{\theta}\varphi(x; \theta) = 0, (x, \theta) \in \Gamma \; \Longrightarrow \; \text{ the collection of vectors } \\ &\left\{ d_{x, \theta} \frac{\partial \varphi(x; \theta)}{\partial \theta_j} : j = 1, \cdots, M \right\} \text{ is linearly independent.}   
\end{align*} 
\item Let $\Lambda$ be a conic $C^{\infty}$ submanifold of $T^{\ast} X \setminus 0$. A nondegenerate phase function $\varphi$ is said to \emph{parametrize} $\Lambda$ if 
\[ \Lambda = \{(x, d_x \varphi(x; \theta)) : d_{\theta} \varphi(x, \theta) = 0, (x, \theta) \in \Gamma \}. \]  
\end{enumerate}
\end{definition}
\begin{definition}\label{def-mp} Let $(\Lambda_0, \Lambda_1)$ be a pair of Lagrangians in $T^{\ast} X \setminus 0$ that 
intersect cleanly in codimension $k$. Let $\lambda_0 \in \Lambda_0 \cap \Lambda_1$ and $\Gamma \subseteq X \times (\R^N\setminus 0)\times\R^k$ an open conic set. A \emph{multiphase function} $\phi$ parametrizing the pair $(\Lambda_0, \Lambda_1)$ is a function $\phi(x;\theta;\sigma)\in C^\infty(\Gamma)$ such that
\begin{enumerate}[(i)]
\item $\phi_0(x;\theta):=\phi(x;\theta;0)$ is a nondegenerate phase
function parametrizing $\Lambda_0$ in a conic neighborhood of $\lambda_0$, and
\item $\phi_1(x;(\theta,\sigma)):=\phi(x;\theta;\sigma)$ is a nondegenerate phase function parametrizing $\Lambda_1$ in a conic neighborhood of $\lambda_0$. \end{enumerate} 
\end{definition}
{\em{Example.}} It is known \cite{meuh,guuh} that any two pairs of cleanly intersecting Lagrangians are microlocally equivalent. One can thus consider the model pair $(\Lambda_0, \Lambda_1)$ in $T^{\ast} \mathbb R^n$ where $\Lambda_0$ and $\Lambda_1$ are the conormal bundles of $\{x = (x_1, \cdots, x_n) = 0\}$ and $\{ x_{k+1} = \cdots = x_n = 0 \}$, resp., so that \begin{align*} \Lambda_0 &= T_0^{\ast} \mathbb R^n = \{(x, \xi) : x = 0 \}, \\ \Lambda_1 &= \{(x, \xi) : x_{k+1} = \cdots = x_n = 0, \xi_1 = \cdots, \xi_k = 0 \}. \end{align*} Then $\varphi(x; \theta', \sigma) = x \cdot (\sigma, \theta')$, with $\theta' \in \mathbb R^{n-k} \setminus 0$, $\sigma \in \mathbb R^k$, is an example of a multiphase function parametrizing $(\Lambda_0, \Lambda_1)$.     In this paper, we will only be concerned with the case $k=1$, i.e., codimension one intersections. 
\begin{definition}\label{def-svs}
The space $S^{\tilde{p},\tilde{l}}\left(X \times (\mathbb R^{N}\setminus 0) \times \mathbb R\right)$ of  \emph{symbol-valued symbols} of orders $\tilde{p},\tilde{l}$,
is the set of functions $a(x;\theta;\sigma)\in C^\infty(X \times (\R^N\setminus \nolinebreak0)\times\R)$ such that, for every relatively compact $K \subseteq X$, non-negative multi-indices $\alpha\in\mathbb Z^N$, $\beta\in\mathbb Z$ and $\gamma\in\mathbb Z^{\dim X}$, the following differential estimates hold:
\begin{equation}\label{ests-svs}
|\partial^\alpha_{\theta}
\partial^\beta_{\sigma}\partial^\gamma_{x}
a(x;\theta;\sigma)| \leq
C_{\alpha,\beta,\gamma,K}\langle\theta,\sigma\rangle^{\tilde{p}-|\alpha|}\langle\sigma\rangle^{\tilde{l}-\beta},
\end{equation}
for all $(x,y)\in K$, where $\langle\xi\rangle=(1+|\xi|^2)^\frac12$, etc., throughout the paper.
\end{definition}
The quantities $\sup_{x, \theta, \sigma} \langle\theta, \sigma\rangle^{|\alpha|-\tilde{p}} \langle\sigma\rangle^{\beta - \tilde{l}} |\partial_{\theta}^{\alpha} \partial_{\sigma}^{\beta} \partial_x^{\gamma}a(x; \theta, \sigma)|$ are referred to as the seminorms for the class $S^{\tilde{p}, \tilde{l}}$. If $|\theta| \geq |\sigma|$ on the support of $a$, then we say that the phase variable $\theta$ is  {\em{dominant}}.   

We  next define the classes of generalized Fourier
integral distributions associated with an intersecting pair of Lagrangians.
\begin{definition}\label{def-ipl-dist} If $\Lambda_0,\Lambda_1\subset
T^*(X \times Y)\setminus 0$ are smooth, conic Lagrangians intersecting cleanly in codimension one, then the space of \emph{generalized} (or \emph{paired Lagrangian}) Fourier integral distributions  of order $p,l\in\R$ associated to
$(\Lambda_0,\Lambda_1)$, denoted $I^{p,l}(X \times Y;\Lambda_0, \Lambda_1)$,
is the set of all locally finite sums of elements of $I^{p+l}(\Lambda_0)+I^p(\Lambda_1)$ and distributions of the form
\begin{equation}\label{eqn-orders}
u(x,y)=\int e^{i
\phi(x,y;\theta;\sigma)} a(x,y;\theta;\sigma) d \sigma d \theta,
\end{equation}
where  $a \in S^{\tilde{p},\tilde{l}}(X\times Y\times (\R^N\setminus
0)\times\R)$, with
\begin{equation}\label{eqn-ipl}
p=\tilde{p}+\tilde{l}+\frac{N+1}{2}-\frac{\dim X + \dim
Y}{4},\quad l=-\tilde{l}-\frac{1}{2},
\end{equation}
and $\phi(x,y;\theta; \sigma)$ is a  multiphase function
parametrizing $(\Lambda_0,
\Lambda_1)$  on a conic neighborhood of a point
$\lambda_0\in\Lambda_0\cap\Lambda_1$.
\end{definition}

{\bf Remark.} Aside from an adjustment in the orders, the paired
Lagrangian spaces are symmetric \cite{meuh,guuh} in $\Lambda_0,\Lambda_1$:
\begin{equation}\label{eqn-sym}
I^{p,l}\left(\Lambda_0,\Lambda_1\right)=
I^{p+l,-l}\left(\Lambda_1,\Lambda_0\right)\supsetneq I^{p+l}(\Lambda_0)+I^p(\Lambda_1).
\end{equation}

Finally, we  define the classes of generalized (or {\it paired
Lagrangian}) Fourier integral operators which are the subject of this paper.
\begin{definition}\label{def-ipl}
(i) If  $C_0,C_1\subset (T^{\ast}X\setminus 0) \times (T^{\ast}Y\setminus 0)$  are smooth, conic canonical relations intersecting cleanly, then  $I^{p,l}(X,Y;C_0,C_1)$  denotes the set of operators $A:{\mathcal E}'(Y)\longrightarrow{ \mathcal D}'(X)$ with Schwartz kernels $K_A(x,y)\in I^{p,l}(X \times Y; C_0',C_1')$.

(ii) In particular, if $C\subset (T^{\ast}Y\setminus 0) \times (T^{\ast}Y\setminus 0$ is a canonical relation intersecting the diagonal
relation $\D$
cleanly, then the members of the associated class $I^{p,l}(\D,C):=I^{p,l}(Y,Y;\D,C)$ are referred to as \emph{pseudodifferential operators with
singular symbols \cite{meuh,guuh}}.
\end{definition}
The spaces $X$ and $Y$ are suppressed if clear from the context.
Equipped with the definitions above, we are now ready to  describe the operators of interest. 

\subsection{Oscillatory  integral representation of operators in $I^{p, l}(\Delta, C_i)$} \label{oscrep} 
The operators we consider will be of the following two forms. 
\subsubsection{The cubic model} \label{cubic-description} For the cubic model, the operator $A$ lies in the class $I^{p, l}(Y, Y; \Delta, C_0)$, where $Y$ is an open bounded subset of $\mathbb R^2$, and $C_0'$ is the conormal bundle of the cubic $\{x_2 - y_2 =(x_1 - y_1)^3 \}$. The kernel of $A \in I^{p,l}(\Delta, C_0)$ can be written as 
\begin{equation}
K_A(x,y) = \int_{\mathbb R^2} e^{i \phi_0(x,y; \theta)} a(x,y; \theta_2; \theta_1)\, d\theta, \label{Skernel0} 
\end{equation}
where the symbol-valued symbol $a \in S^{p+\frac{1}{2}, l-\frac{1}{2}}(Y \times Y \times (\mathbb R_{\theta_2} \setminus 0) \times \mathbb R_{\theta_1})$, and the multi-phase function $\phi_0$ is given by  
 \begin{equation}\label{phi_0}
\phi_0(x,y,\theta)=(x_1-y_1)\theta_1+(x_2-y_2-(x_1-y_1)^3)\theta_2. 
\end{equation}
The reader may verify that $\phi_0(x,y;\theta_2;\theta_1)$ is a multi-phase function parametrizing, as in Def. \ref{def-mp}, the cleanly intersecting pair of canonical relations $(C_0, \Delta)$, in that order. The dominant phase variable is $\theta_2$.  

\subsubsection{Normal operators for seismic imaging}\label{noropseis}
We now turn to the operator classes $I^{p,l}(Y, Y; \Delta, C_j)\, j=1,2$, where $(\D,C_j)$ are the geometries  that arise in the analysis of the normal operators $F^{\ast}F$  in seismic imaging, as described in the introduction. Here $Y$ is a bounded open subset of $\mathbb R^n$, $n \geq 3$ ($n=3$ for the physical problem).

For the single source problem, the Melrose-Taylor 
normal form  \cite{meta} for   folding canonical relations, such as the one to which $F$ is associated, 
 gives rise to a particularly simple normal form for the pair $(\D,C_1)$,
namely $\left(\Delta, C_0 \times \Delta_{T^*\R^{n-2}}\right)$.
This  allows for the estimates for $I^{p,l}(\D,C_1)$ to be proven as for $I^{p,l}(\D,C_0)$; see \S\S\ref{subsec-remark}.

In contrast, the
marine data geometry is complicated by the fact that there is only an
\emph{approximate} normal form for folded cross caps. It was shown
in \cite{FG} that, microlocally near any point
in $\Delta \cap C_2$  the pair $(\Delta, C_2)$
can be parametrized, in that order, by a multi-phase function of the form
\begin{align}\label{mg-phi}
\phi(x,y;\xi;\rho)&=(x-y)\cdot\xi + \frac{\rho}{\xi_1}\left( \xi_n -
\frac{(x_n+y_n)^2}4 \xi_1 -P\right), \text{ on } \\ 
\big\{\,|\xi_1| &\geq \frac{1}{2}(|\xi| + |\rho|)\,\big\}. \label{xi1-dominant} 
\end{align}
Here, $P = P(x_n,y,\xi')$ is an unknown smooth function  on
$\R\times\R^n\times(\R^{n-1}\setminus 0)$, homogeneous of degree one
in $\xi'= (\xi_1, \cdots, \xi_{n-1})$ and satisfying
\begin{equation}\label{mg-P}
P|_{\{x_n=y_n\}}=0,\quad \nabla P|_{\{x_n=y_n\}}=0
\end{equation}
and
\begin{equation}\label{mg-P2}
||\nabla^2_{x_n,y_n} P|| \ll |\xi_1|
\end{equation}
microlocally, i.e., one can arrange for the left hand side of (\ref{mg-P2}) to vanish at any chosen base point, and thus be small nearby.  Here, $\xi_1$ is a dominant phase variable. Note that we have translated the formulas in \cite{FG}, which were written in $n-1$ dimensions, to an $n$-dimensional setting.

However, in order to handle this geometry in a way similar to the cubic model $(\D,C_0)$, we want a multiphase function that parametrizes the pair $(\Delta,C_2)$  in the reverse order. A prescription from \cite{men}, which we now describe, allows us to compute such a phase function explicitly. We observe that $\Delta'$ is the flowout from $C_2'\cap\Delta'$ by the Hamiltonian  vector field $H_q$ on $T^*\left(\R^n\times\R^n\right)$, where $q(x,\xi,y,\eta)=(x_n-y_n)\xi_1$. Starting with $\psi_0(x,y;\xi,\rho)=\phi$ as in
(\ref{mg-phi}), which parametrizes $C_2$, we solve for $\tilde\psi(x,y;\xi,\rho;s)$ satisfying
\begin{equation*}
\frac{\p\tilde\psi}{\p s} = (x_n-y_n)\frac{\p \tilde\psi}{\p 
x_1},\quad \tilde\psi|_{s=0}=\psi_0,
\end{equation*}
and find that
\[\tilde\psi=(x-y)\cdot\xi+\frac{\xi_n\rho}{\xi_1} -\frac14(x_n+y_n)^2\rho -
\frac{\rho}{\xi_1} P +s(x_n-y_n)\xi_1. \]
Finally, letting $s=\frac{\sigma}{\xi_1}$, this yields the multi-phase function
\begin{equation}\label{mg-psi}
\psi(x,y;(\xi,\rho);\sigma)=(x-y)\cdot \xi +\frac{\xi_n\rho}{\xi_1} -
\frac14 (x_n+y_n)^2\rho -\frac{\rho}{\xi_1} P +(x_n-y_n)\sigma,
\end{equation}
which parametrizes $(C_2,\Delta)$ in that order according to Def. \ref{def-mp}. By Def. \ref{def-ipl-dist}, and the relations (\ref{eqn-sym}) and (\ref{eqn-ipl}), the kernel of an operator
$A\in I^{p,l}(\Delta,C_2)= I^{p+l,-l}(C_2,\Delta)$ has a representation, modulo $I^{p+l}(\Delta) + I^p(C_2)$  of the form 
\begin{equation}\label{mg-KA}
K_A(x,y)=\int e^{i\psi(x,y;\xi,\rho;\sigma)}  a(x,y;\xi,\rho;\sigma) 
d\xi d\rho d\sigma,
\end{equation}
with $a\in S^{p-\frac12,l-\frac12}\left(Y \times Y \times 
(\R^{n+1}_{\xi,\rho}\setminus 0)\times\R_\sigma\right)$, and supp$(a)
\subseteq \{ |\sigma| \leq |\langle \xi, \rho \rangle| \}$.

\subsection{A technical tool: the method of stationary phase}
In the subsequent sections, especially \S \ref{sec-gendecomp} and \S
\ref{sec-norm}, we will need to estimate repeatedly the oscillatory
integrals that arise from (\ref{Skernel0}) and (\ref{mg-KA}) as the Schwartz kernels of the compositions $B^{\ast}A$
and $BA^{\ast}$, where $A, B \in I^{p,l}(\Delta, C_j)$. We will use the method of
stationary phase to describe the asymptotic behavior of integrals
of the form 
\[ I(a, \lambda) = \int_{\mathbb R^n} e^{i \lambda \varphi(w, a)} g(w,
a, \lambda) \, dw, \qquad a \in \mathbb R^p \]
as $\lambda \rightarrow \infty$; for the sake of completeness, we state it in the form needed. The phase function $\varphi$ is
assumed to be real-valued and smooth, while the amplitude function $g$
is assumed to be smooth of compact support in the variables $(w,a) \in
\mathbb R^n \times \mathbb R^p$ with the following growth restriction: there exists some $\eta < \frac{1}{2}$ such that 
\begin{equation} \left(\frac{\partial}{\partial w}\right)^{\alpha} g = O(\lambda^{m
+ \eta |\alpha|}) \;  \text{ as
$\lambda \rightarrow \infty$, uniformly in $(w,a)$}. \label{g-est} \end{equation} We state
here the principle of stationary phase in a form that we need. 
\begin{lemma}[\cite{du}, p.14]
Let $\varphi(w, a) = \frac{1}{2}\langle Q(a)w, w \rangle$, where $Q$
is a real $n \times n$ nonsingular, symmetric matrix depending
continuously on $a$, and let $g$ be an amplitude satisfying
(\ref{g-est}) above. Then $I(a, \lambda)$ has an asymptotic expansion 
\[ I(a, \lambda) \sim \left( \frac{2 \pi}{\lambda}\right)^{\frac{n}{2}}
\bigl|\det (Q(a)) \bigr|^{-\frac{1}{2}} e^{\frac{\pi i}{4}
\text{sgn}(Q(a))} \sum_{r=0}^{\infty} \frac{(R^r g)(0, a,
\lambda)}{r!} \lambda^{-r} \]   
as $\lambda \rightarrow \infty$, uniformly in $a$. Here 
\[ R = \frac{i}{2} \left \langle Q(a)^{-1} \frac{\partial}{\partial w},
\frac{\partial}{\partial w} \right \rangle \] 
which is a second-order partial differential operator in $w$. 
\label{SP-lemma}\end{lemma}
In the sequel, justifying an application of stationary phase will
involve identifying an oscillatory integral with (quadratic) phase $\lambda
\varphi$ and amplitude $g$, and verifying that (a) the Hessian of
$\lambda \varphi$ evaluated at its critical point is nonsingular, and
(b) the quantity $(\lambda^{-1}R)^r g$ evaluated at this critical
point decays exponentially in $r$. Also, while the notation $\sim$ was used in the above lemma to denote an asymptotic expansion of $I(a,\lambda)$, in subsequent applications we will use $\sim$ to represent just the leading order term of the asymptotic expansion, with the understanding that the latter also provides the sharp size estimate for $I(a, \lambda)$ for $\lambda \gg 1$.   

\section{Sharpness of Theorem \ref{thm-main}}\label{sec-sharp}
We start by showing that, at least for $I^{p,l}(\D,C_0)$,  the estimates in Thm.\ \ref{thm-main} cannot be improved in general. For $l < - \frac{1}{2}$, the inclusion $I^p(C_0) \subset I^{p,l}(\Delta, C_0)$, cf. (\ref{eqn-sym}),  implies that the Sobolev mapping index $r$ of $I^{p,l}(\Delta, C_0)$ must be at least as large as that of (\ref{est meta}) for $I^p(C_0)$, i.e., $r \geq p + \frac{1}{6}$. On the other hand, the fact hat $\Psi^{p+l}(\R^2)= I^{p+l}(\Delta)\subset I^{p,l}(\Delta, C_0)$ implies that we must have $r\ge p+l$.

In the critical interval $-\frac12 \le l <\frac12$, we demonstrate the optimality of 
the estimates of Thm.\ \ref{thm-main} for the fractional integral operators along $(t,t^3)$. For $-\frac12 < l<\frac12$,  define, as in (\ref{def-frac}),
\be\label{def-tl}
\mathcal J_{l+\frac12}f(x)=f*\left(\delta(x_2-x_1^3)\cdot\frac{\chi(x_1)\,
dx_1}{|x_1|^{l+\frac12}}\right),
\end{equation}
where $\chi\in\coi(\R)$, $\chi(t)\equiv 1$ for $|t|\le\frac12$. Writing the kernel of the operator as an oscillatory integral, it follows from the discussion in \S\S\ref{cubic-description} that $\mathcal J_{l+\frac12}\in I^{-\frac12,l}(\D,C_0)$.  Thm.\ \ref{thm-main} therefore implies that
$\mathcal J_{l+\frac12}:H^s_{comp}(\R^2)\lra H^{s-r}_{loc}(\R^2)$ for
$r=\frac{2l-1}6$. That this cannot be improved is seen as follows.
Since the Fourier multiplier of $\mathcal J_{l+\frac12}$ is
\be
m_l(\xi)=\int_\R e^{-i(\xi_1 t + \xi_2 t^3)} \frac{\chi(t)\,
dt}{|t|^{l+\frac12}},
\end{equation}
it suffices to show that \[ \sup_{\xi\in\R^2} |m(\xi)|
|\xi|^{\frac{1-2l}6}\ge c>0. \] 
Using the substitution $u=\xi_2 t^3,\, du=(3\xi_2^{1/3}
u^{2/3})^{-1}dt$, one has
\begin{equation} \begin{aligned} m_l(0,\xi_2)&= \int e^{-i\xi_2 t^3} \frac{\chi(t)\, dt}{|t|^{l+\frac12}}\\
&= c |\xi_2|^{\frac{2l-1}6}\int e^{-iu} \frac{\chi\left(\xi_2^{-\frac13}
u^{\frac13}\right)\, du}{|u|^{\frac{2l+5}6}}. \end{aligned} \label{int-last}
\end{equation} 
Since the improper integral $\int e^{-iu} |u|^{-\frac{2l+5}6} \, du$ converges to a nonzero value, the same is true for the  integral in (\ref{int-last}).

For $l=-\frac12$, we replace $|x_1|^{-0}$ by $\log |x_1|$ in (\ref{def-tl}), which still defines an operator $\mathcal J_{0}\in I^{-\frac12,-\frac12}(\D,\tC_0)$
since $\log |x_1|$ is conormal of order $-1$ for $\{x_1=0\}$ on $\R$.  The
corresponding logarithmic divergence of
$|m_{-\frac12}(0,\xi_2)|\cdot |\xi_2|^{\frac13}$ then implies that
one cannot in general eliminate the loss of $\e$ derivatives when $l=-\frac12$ in Thm. \ref{thm-main}.

\section{The canonical graph case}\label{sec-graph}
As a warmup for the proof of Thm.\ \ref{thm-main}, and to provide a background for the decompositions that are needed, we first prove an analogous result under the optimal, nondegenerate microlocal geometry, namely when $C$ is a canonical graph. While the intersection $C\cap\Delta$ can be quite singular, we restrict ourselves to the configuration closest to that of Thm. \ref{thm-main}, namely when $C$ intersects $\D$ cleanly, say in codimension $k$. For example, if $\gamma$ in (\ref{eqn-hil}) is replaced with the parabola $(t,t^2)$, then  that $\mathcal H\in I^{-\frac12,\frac12}(\D, C)$ for a local canonical graph $C$, although  for families of curves in higher dimensions, the canonical graph condition must always be violated somewhere. We mention in passing that singular Radon transforms associated with a family of $k$-dimensional surfaces in $\R^n$ belong to the class $I^{-\frac{k}2,\frac{k}2}(\Delta, C)$, where $C$ is a local canonical graph under the rotational curvature condition \cite{PS1}, and with $\Delta$ intersecting $C$ cleanly in codimension $k$ \cite{gruh2}. Estimates for such operators are in \cite{NSW,GeSt,PS1,gruh2,PS2,seW03}; see also \cite{CuThesis} for worst-case estimates in degenerate cases. 

The following result is similar to Thm. \ref{thm-main}, but much easier to prove due to the nondegeneracy of $C$.
For $C$ a conormal bundle, related estimates were found by Seeger and Wainger \cite{seW03}.

\begin{theorem}\label{thm-cg}
Let $Y$ be a manifold of dimension $n$. Let $C\subset (T^*Y\setminus 0)\times (T^*Y\setminus 0)$ be a local canonical graph  intersecting $\D$ cleanly in codimension  $k,\quad1\le k\le n-1$. If $A\in
I^{p,l}(\D, C)$, then $A:H^s_{comp}\lra H^{s-r}_{loc}$ for
\begin{equation}\label{est-cg}
r=\begin{cases} p, & l<-\frac{k}2,\\
p+\epsilon, & l=-\frac{k}2,\, \forall \e>0,\\
p+\frac{2l+k}4, & -\frac{k}2<l<\frac{k}2,\\
p+l, & l\ge \frac{k}2.
\end{cases}
\end{equation}
\end{theorem}

\begin{proof} We  prove the theorem using a standard parabolic cutoff.
Parabolic cutoffs  in various guises  have been used often in microlocal analysis; the argument  given here is a direct adaptation of that in
\cite{gruh2}, but related ideas are also in \cite{NSW,PS2,mel}. Since $A\in I^{p,l}(\Delta,C)= I^{p+l,-l}(C,\Delta)$ by (\ref{eqn-sym}), we can assume that $A$ has  kernel

\begin{equation} K_A(x,y)=\int_{\R^{N+1}} e^{i\phi(x,y;\theta;\sigma)}
a(x,y;\theta;\sigma) d\theta d\sigma, \label{Skernel-cg} \end{equation}
where $\phi$ is a multiphase function parametrizing the cleanly
intersecting pair $(C,\Delta)$ and $a\in S^{p-\frac{N-n}2, l-\frac{k}2} (Y\times Y\times (\R^N\setminus 0)\times\R)$. Decompose

$$a=\chi\left(\frac{\la\sigma\ra}{\la\theta\ra^\frac12}\right)\cdot a +
\left(1-\chi\left(\frac{\la\sigma\ra}{\la\theta\ra^\frac12}\right)\right)\cdot
a := a_0+a_1,$$
with $A=A_0+A_1$ the corresponding decomposition of $A$.
Since $\la\sigma\ra\ge c \la\theta\ra^\frac12$ on
$supp(a_1)$, the differential estimates on $a_1 = a_1(x, y; (\theta, \sigma))$ ensure that 
it is a standard (non product-type) symbol,
\begin{equation}\label{a1}
a_1\in
\begin{cases} S^{p-\frac{N-n}2+\frac12(l-\frac{k}2)}_{\frac12,0}\left(Y \times Y \times (\mathbb R^{N+k} \setminus 0)\right),\quad &
l<\frac12,\\ S^{p+l-\frac{N-n+k}2}_{\frac12,0}\left(Y \times Y \times (\mathbb R^{N+k} \setminus 0)\right),\quad & l\ge\frac12.
\end{cases}
\end{equation}
Since $\phi(x,y;(\theta,\sigma))$ parametrizes $\Delta$, this implies that

\begin{equation}\label{A1}
A_1\in \Psi_{\frac12,\frac12}^{\max(p+\frac{2l+k}4,p+l)}(Y). 
\end{equation}

We now turn to $A_0$, whose Schwartz kernel is of the form (\ref{Skernel-cg}), but with the amplitude $a$ replaced by $a_0$. Recalling that $\phi(x,y; \theta, 0)$ parametrizes $C$ and reasoning as in \cite[Prop. 2.1]{gruh2}, let us integrate out $\sigma$ and write 

\begin{align*} K_{A_0}(x, y) &= \int_{\mathbb R^N} e^{i \phi(x,y; \theta, 0)} b_0(x,y; \theta) d\theta, \quad \text{where} \\ 
b_0(x,y; \theta) &= \int e^{i (\phi(x,y; \theta, \sigma) - \phi(x,y; \theta, 0))} a_0(x, y; \theta, \sigma) \, d\sigma. 
\end{align*} 
Keeping in mind the symbol estimates for $a_0$ and the homogeneity of $\phi$ in $(\theta, \sigma)$, the argument in \cite{gruh2} carries over with minor modifications to yield
\begin{equation} b_0 \in S^{p - \frac{N-n}{2} + \frac{1}{2}(l + \frac{k}{2})_{\tilde{+}}}_{\frac{1}{2}, \frac{1}{2}}, \quad \text{ and hence } \quad A_0\in I_{\frac12,\frac12}^{p+\frac12 (l+\frac{k}2)_{\tilde{+}}}(C). \label{A0}\end{equation} 
Here, and throughout, for $t\in\R$, we set
\begin{equation} t_{\tilde{+}}= \begin{cases} t  &\text{ if } t>0,  \\  0 &\text{ if }  t<0, \text{ and } \\  \text{ any } \epsilon>0 & \text{ for }  t=0. \end{cases} \label{plus-tilde} \end{equation} 
Combining (\ref{A0}), (\ref{A1}) it follows from the Calder\'on-Vaillancourt
theorem and its extension to FIOs of type $(\frac12,\frac12)$ associated with
local canonical graphs \cite{BEALS,gruh2} that
$A:H^s_{comp}\lra H^{s-r}_{loc}$ for
$$r\ge\max\left(p+\frac12 (l+\frac{k}2)_{\tilde{+}},p+\frac{2l+k}4,p+l\right),$$
which yields  (\ref{est-cg}). 
The calculations of \S \ref{sec-sharp} are easily modified for the curve $(t,t^2)$ to show that Theorem \ref{thm-cg} is  sharp, at least for $k=1$.
\end{proof}
When $C$ is a folding canonical relation,  more involved decompositions of $A$, requiring additional knowledge of the multiphase functions,  are needed to prove the analogous estimates. \S\ref{sec-gendecomp},\S\ref{sec-norm} and \S\ref{sec-cubic}  are devoted to dealing with  the cases when $C$ has one of the forms described in \S\S\ref{oscrep}.

\section{Mixed Parabolic-Phong-Stein Decompositions for
$I^{p,l}(\Delta, C_0)$} 
\label{sec-gendecomp}
In the previous section, we  described and applied  the method of
parabolic cutoffs in the context of FIOs associated with canonical
graphs. There is also a method introduced by Phong and Stein \cite{PS3,PS4} for proving $L^2$ estimates using dyadic decompositions in phase space; 
this, along with
variations, has proven to be a powerful technique for controlling oscillatory integral
operators and FIOs with degenerate canonical relations
\cite{se93,Cu,Co,grse2}. In the present context, our analysis of  operators in
$I^{p,l}(\Delta, C_j)$, $j = 0,2$, will combine   these two decomposition
strategies. In this section, we describe this  for the class
$I^{p,l}(\Delta, C_0)$, where we recall that $C_0 = N^{\ast}\{x_2 -
y_2 = (x_1 - y_1)^3 \}'$. Appropriate modifications of these
decompositions for  $I^{p,l}(\Delta, C_2)$ will be dealt
with in \S \ref{sec-norm}.   

\subsection{Components of the operator} \label{gen-decomp} Let $\{\chi_j : j \geq 0 \}$ be a non-homogeneous dyadic partition of unity on $\mathbb R$ such that 
\begin{equation} \begin{aligned} \supp(\chi_0) &\subset \{|t|\le 2\}, \\ \supp(\chi_j) &\subset\{2^{j-1}\le |t|\le
2^{j+1}\}, \quad j\ge 1, \text{ and } \\    |\chi_j^{(m)}| &\le C_m
2^{-mj},\quad \text{ for all } m \geq 1. \end{aligned} \label{pou-chi} \end{equation} 
Setting $\chi_{jk}(\theta)=\chi_k(\theta_1)\chi_j(\theta_2)$, we note that
$\{\chi_{jk}| 0\le k\le j<\infty\}$ is a bounded family in $S^{0,0}$
with respect to the seminorms of (\ref{ests-svs}).  Furthermore, we set
\[\chi_\infty:=1-\sum_{j \geq 0}\sum_{k=0}^j\chi_{jk}, \] so that
$\chi_{\infty} \in
S^0_{1,0}$ and is supported in $\{|\theta_1|\ge \frac{1}{2} |\theta_2|\}$.
Given any $A \in I^{p,l}(\Delta, C_0)$ whose Schwartz kernel is of the
form (\ref{Skernel0}), its amplitude $a(x,y;\theta_2; \theta_1)$ decomposes as
\begin{align*}
a(x,y; \theta_2, \theta_1) &= \left(\sum_{j \geq 0} \sum_{k=0}^j \chi_{jk}(\theta) +\chi_\infty (\theta) \right)\cdot a(x,y; \theta_2; \theta_1) \\
&:= \sum_{j=0}^{\infty}\sum_{k=0}^{j} a_{jk}(x,y; \theta_2; \theta_1) +a_\infty(x,y; \theta_2; \theta_1). 
\end{align*}
Thus $\{a_{jk}\}$ is a bounded family in $S^{p+\frac12,l-\frac12}$, with 
\[ \supp(a_{jk}) \subseteq \{(x, y; \theta) : 2^{j-1} \leq |\theta_2| \leq 2^{j+1},\; 2^{k-1} \leq |\theta_1| \leq 2^{k+1}  \}, \]
and $a_\infty\in S^{p+l}_{1,0}$.
The decomposition of the amplitude $a$ induces a decomposition of the operator $A$,
$$A =\sum_{j=0}^\infty\sum_{0\le k\le j} A_{jk} +A_\infty,$$
where the $A_{jk}$ and  $A_{\infty}$  are operators whose Schwartz kernels and multi-phase functions are of the form described in  \S\S\ref{cubic-description} (see (\ref{Skernel0}) and (\ref{phi_0})), but with the amplitudes replaced by $a_{jk}$ and $a_{\infty}$, resp. Thus, $A_\infty\in I^{p+l}(\Delta_{T^*\R^2})=\Psi^{p+l}(\R^2)$, 
and $\left\{ A_{jk} : 0\le k\le j<\infty\right\}\subseteq I^{p,l}(\D,C_0)$. It follows from standard $\Psi$DO estimates that $A_{\infty} : H^s_{\text{comp}} \rightarrow H^{s - (p+l)}_{\text{loc}}$, which is at least as smoothing as what is being claimed for the entire operator $A$ in Theorem \ref{thm-main}. We therefore disregard $A_{\infty}$ in the sequel. 
\medskip

For a fixed constant $\delta\in(0,1)$ to be set later (we will eventually take $\delta\in [\frac13,\frac12)$), we group the summands in $A - A_{\infty}$ as follows,   
\begin{align} A - A_{\infty} &= A_0 +
\sum_{j=1}^{\infty} A_j, \text{ where } \label{comp1} \\
A_0 = \sum_{j=0}^{\infty} \sum_{k=0}^{[\delta j]} A_{jk}, &
\quad  A_j = \sum_{k = [\delta j]+1}^{j}A_{jk}, \; j \geq 1. \label{comp2}
\end{align}
In this section we identify $A_0$ as a standard FIO
associated with the  folding canonical relation $C_0$ and
record several almost orthogonality properties involving the $A_j$,
$A_{j\infty}$ and $A_{jk}$. These facts will be used heavily
in the proof of $L^2$-Sobolev estimates.


\begin{figure}
\hskip-2.75cm
\includegraphics[width=1.2\linewidth]{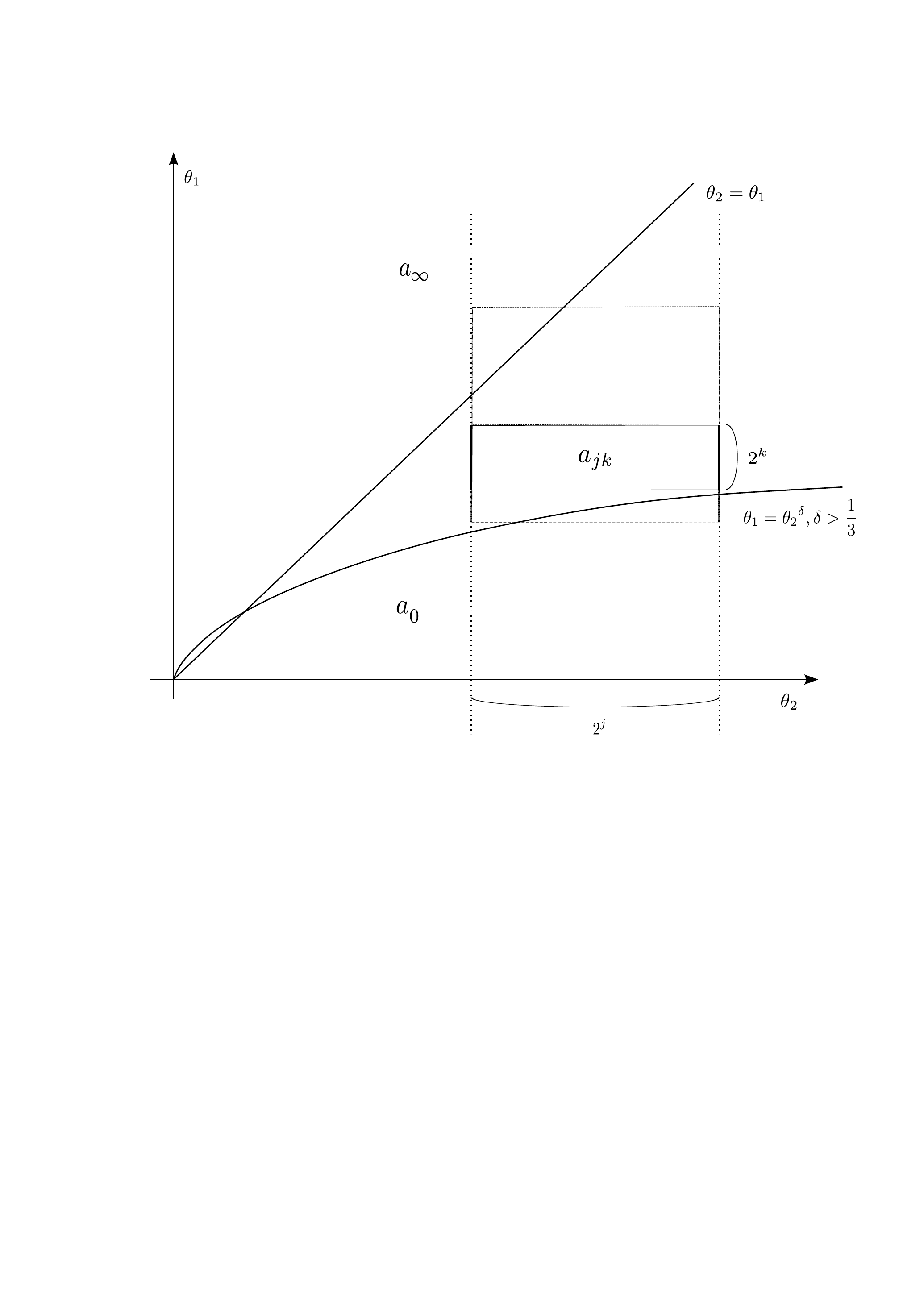}
\vskip-9cm
\caption{Decomposition of $(\theta_2,\theta_1)$ space.}
\end{figure} 




\subsection{Bounds for $A_0$}
\begin{lemma}
For any $0 < \delta < \frac{1}{2}$, the operator $A_0$ lies in the class $I^m_{1 - \delta, \delta}(C_0)$, where $m = p + \delta(l + \frac{1}{2})_{\tilde{+}}$, and hence maps $H^{s}_{\text{comp}}(\mathbb R^2)$ boundedly into $H^{s - r_0}_{\text{loc}}(\mathbb R^2)$, where 
\begin{equation} \label{r0} r_0 \geq p + \frac{1}{6} + \delta \bigl(l + \frac{1}{2} \bigr)_{\tilde{+}}. \end{equation}
Here $t_{\tilde{+}}$ is as   in (\ref{plus-tilde}).  \label{boundsA0}
\end{lemma} 
\begin{proof}
The kernel  of $A_0$ is of the form (\ref{Skernel0}), but with the amplitude replaced by 
\[a_0 = \sum_{j=0}^{\infty} \sum_{k=0}^{[\delta j]} a_{jk}. \]
We observe that $\supp(a_0) \subseteq \{\langle \theta_1 \rangle \leq C \langle \theta_2 \rangle^{\delta} \}$, and that $a_0$ satisfies the same differential inequalities as $\tilde{\chi}\left({\la\theta_1\ra}{\la\theta_2\ra^{-\delta}}\right)\cdot a(x,y;\theta_2;\theta_1)$, for some $\tilde{\chi} \in C_0^{\infty}(\mathbb R)$. 
Integrating  out $\theta_1$ in the oscillatory representation of
$K_{A_0}$ as in the proof of Theorem \ref{thm-cg}, we see that 
\begin{align*} K_{A_0}(x,y) &= \int e^{i \theta_2(x_2 - y_2 - (x_1 - y_1)^3)} \tilde{a}_0(x, y; \theta_2) \, d\theta_2, \quad \text{ where } \\ \tilde{a}_0(x, y; \theta_2) &= \int e^{i(x_1 - y_1) \theta_1} a_0(x, y ; \theta_2, \theta_1) \, d\theta_1.   \end{align*} 
It is easily verified from the support and differentiability properties of $a_0$, using arguments similar to those used in Theorem \ref{thm-cg} that $\tilde{a}_0 \in S_{1, \delta}^{m+\frac{1}{2}}$, where $m$ is as in the statement of the lemma. In other words $A_0$ is an FIO associated with the two-sided folding canonical relation $C_0$, of
order $m +\frac12 +\frac12-\frac44=m$ and  with an amplitude of type
$(1, \delta)$. Since useful FIO classes need to be invariant under changes of variables and phase
functions, we may as well consider (cf. \cite{gruh2}) the type as being
$(1-\delta,\delta)$.   Now, $1-\delta>\frac12$ since 
$\delta<\frac12$,  so (\ref{r0}) follows from \cite{meta}.
\end{proof}

\subsection{Almost orthogonality}
The next two lemmas deal with the almost orthogonality among the components $\{A_{jk} \}$. 
Both for clarity and for possible future use, we will treat the $A^*$ in $A^*A$ as the adjoint of a general operator $B\in I^{p',l'}(\D, C_0)$.

\begin{lemma}[Almost orthogonality in $j$]
Fix $0 < \delta < 1$, and orders $p, p', l, l' \in \mathbb
R$. Given $A \in I^{p,l}(\Delta, C_0)$ and $B \in I^{p',l'}(\Delta,C_0)$, we decompose $A$ and $B$ as described in $\S\S$\ref{gen-decomp} using the same value of $\delta$ for both. Let $\{A_{jk} : 0 \leq k \leq j \}$ and $\{B_{jk} : 0 \leq k \leq j \}$ denote the components of $A$ and $B$ respectively (see equations (\ref{comp1}) and (\ref{comp2})). Then for every $N \geq 1$, there exists a constant $C_N > 0$ such that for all $j, j' \geq 1$, $|j - j'| \geq 3$, $\delta j \leq k \leq j$, $\delta j' \leq k' \leq j'$, the following estimate holds: 
\begin{equation} \bigl| \bigl| B_{j'k'}^{\ast} A_{jk} \bigr| \bigr|_{L^2 \rightarrow L^2} \leq C_N 2^{-N \max(j,j')}. \label{aorth1-ineq1} \end{equation}
In particular, if $|j-j'| \geq 3$, $j, j' \geq 1$, then 
\begin{equation} ||B_{j'}^{\ast} A_j||_{L^2 \rightarrow L^2} \leq C_N
2^{-N \max(j,j')}. \label{aorth1-ineq2} \end{equation} \label{aorth1}
The same estimates also hold for $B_{j'k'} A_{jk}^{\ast}$ and
$B_{j'}A_{j}^{\ast}$.  
\end{lemma}

\begin{proof}
We only give the proof for $B^{\ast}A$, that for $BA^{\ast}$ being
similar. The Schwartz kernel of the composition $B_{j'k'}^{\ast} A_{jk}$ has the form 
\begin{align}\label{Kjpj}
K_{B_{j'k'}^* A_{jk}}(x,y)=\int e^{i\Phi(x,y;z;\theta;\theta')}\overline{b}_{j'k'}(z,x;\theta') a_{jk}(z,y;\theta) d\theta d\theta' dz,
\end{align}
where the phase function $\Phi$ is given by 
\begin{multline}
\Phi(x, y, z; \theta, \theta')  =  (z_1-y_1)\theta_1 - (z_1-x_1)\theta_1'\\
 + (z_2-y_2-(z_1-y_1)^3)\theta_2 - (z_2-x_2- (z_1-x_1)^3)\theta_2'.
 \label{Phi} \end{multline} 
We note that $d_{z_2}\Phi = \theta_2-\theta_2'$, so that for $|j - j'| \geq 3$ and on $\supp(\overline{b}_{j'k'}\cdot a_{jk})$, \[|d_{z_2}\Phi|\geq (2^{\max(j,j')-1} - 2^{\min(j,j') + 1}) \geq \frac{1}{4}2^{\max(j,j')}. \] Integrating by parts
$N$ times in $z_2$, integrating in all
variables, and applying the support and differentiability
properties of $a_{jk}$ and $b_{j'k'}$ throughout we obtain
\begin{multline}\label{est-jjp}
\left|K_{B_{j'k'}^* A_{jk}}(x,y)\right| \leq C_N
2^{(j+k)+(j' + k')} \times \\ 2^{(p+\frac12)j+(p' + \frac{1}{2})j'
+(l - \frac{1}{2})k+ (l' - \frac{1}{2})k'-N\max(j,j')}.
\end{multline}
Integrating this over the compact support in either $x$ or $y$ yields the estimate 
\[ \sup_x \int \left|K_{B_{j'k'}^{\ast} A_{jk}}(x,y)\right| dy + \sup_y \int \left|K_{B_{j'k'}^{\ast} A_{jk}}(x,y)\right| dx  \leq 
C_N 2^{(M-N)\max(j,j')}, \] where $M = 4 + |p + \frac{1}{2}| + |p' +
\frac{1}{2}| + |l - \frac{1}{2}| + |l' - \frac{1}{2}|$ is a fixed
constant. Schur's lemma then implies the same bound for $||B^*_{j'k'}
A_{jk}||_{L^2 \rightarrow L^2}$. Since $N$ is arbitrary, the proof of (\ref{aorth1-ineq1})
is complete. The inequality (\ref{aorth1-ineq2}) follows from (\ref{aorth1-ineq1}) simply by summing in $k$ and $k'$. 
\end{proof}

\begin{lemma}[Almost orthogonality in $k$]
Under the same hypotheses as Lemma \ref{aorth1}, the following conclusions hold. 
\begin{enumerate}[(a)]
\item There exists a constant $C > 0$ such that for all $j \geq 1$ and
$[\delta j] + 1 \leq k \leq j$, \begin{equation} ||A_{jk}||_{L^2
\rightarrow L^2}  \leq C 2^{(p + \frac{1}{2})j + (l - \frac{1}{2})k}. \label{AjkL2} \end{equation} A similar statement holds for $B_{jk}$. \label{parta}
\item Now suppose $\delta \geq \frac{1}{3}$. Then, for any $N \geq 1$, there
exists a constant $C_N > 0$ such that for all $j, j', k, k' \geq 1$,
with $|j-j'| \leq 2$, $|k - k'| \geq 3$, $[\delta j] + 1 \leq k \leq
j$ and $[\delta j'] + 1 \leq k' \leq j'$, the following estimate holds:  
\begin{equation} ||B_{j'k'}^{\ast} A_{jk}||_{L^2 \rightarrow L^2} \leq
C_N 2^{(p + p'+ 1)j + (l - \frac{1}{2})k + (l' - \frac{1}{2})k'}
2^{-|k - k'|N}. \label{aok2} \end{equation} The same estimate holds
for $B_{j'k'}A_{jk}^{\ast}$.  
\label{partb}
\end{enumerate} 
\label{aok}
\end{lemma}

\begin{proof}
For any choice of indices $j, j', k, k'$ with $[\delta j]+1 \leq k
\leq j$ and $[\delta j'] + 1 \leq k' \leq j'$, we make a preliminary
simplification of the integral kernel $K_{B_{j'k'}^{\ast}A_{jk}}$ that
will be useful in the proof of both (\ref{parta}) and
(\ref{partb}). Going back to the oscillatory representation
(\ref{Kjpj}) with $\Phi$ as in (\ref{Phi}), we perform stationary
phase in $(\theta_2', z_2)$ on this oscillatory integral. This step is
justified by Lemma \ref{SP-lemma}, with $w = (\theta_2', z_2)$, $a =
(z_1, x, y, \theta, \theta_1')$ and $\lambda \varphi = \Phi$ in the
notation of that lemma, so that $(\lambda^{-1}R)^r[\bar{b}_{j'k'}
a_{jk}] = (d_{\theta_2' z_2})^r[\bar{b}_{j'k'} a_{jk}] = O(2^{-jr}) \ll 1$.
We are thus led to the estimate  
{\allowdisplaybreaks \begin{multline} \label{eqn-jkpk}
 \left|K_{B_{j'k'}^{\ast}A_{jk}}(x,y) \right| \sim \left|\int e^{i\tP(x,y;z_1;\theta;\tho')} \overline{b}_{j'k'}(\tilde{z},x;\tht;\tho') a_{jk}(\tilde{z},y;\theta)
dz_1 d\tho' d\theta \right|,
\end{multline}}
where $\tilde{z}=(x_1, x_2+(z_1-x_1)^3)$ and
\be\label{eqn-Phit}
\tP=(z_1-y_1)\tho-(z_1-x_1)\tho'+(x_2-y_2+(z_1-x_1)^3-(z_1-y_1)^3)\tht.
\ee
Several integration by parts will now be performed on the integral
in (\ref{eqn-jkpk}). In the subsequent computation, $N \geq 1$ will denote an
arbitrarily large integer and $u_{jj'kk'}(x,y, z_1; \theta_2, \theta_1,
\theta_1')$ any function satisfying similar size and
differentiability estimates as $\bar{b}_{j'k'}(\tilde{z}, x; \theta_2,
\theta_1') \cdot a_{jk}(\tilde{z}, y; \theta)$, possibly with different implicit constants. The value of $N$ or the exact functional form of $u_{jj'kk'}$ may vary from one occurrence to the next.

For part (\ref{parta}), we set $j=j'$, $k=k'$, integrate by parts $2N$ times in (\ref{eqn-jkpk}) with respect to each of the variables $\theta_2$, $\theta_1$, $\theta_1'$, finally obtaining
\begin{multline*} \bigl| K_{B_{jk}^{\ast}A_{jk}}(x,y) \bigr| \le C_N \int \left[ 1 + (2^k(z_1 - y_1))^2\right]^{-N} \left[ 1 + (2^k(z_1 - x_1))^2\right]^{-N} \\ \times \bigl[1 + \bigl(2^j(x_2 - y_2 +(z_1 - x_1)^3 - (z_1 - y_1)^3) \bigr)^2 \bigr]^{-N} u_{jjkk} \, dz_1 \, d\theta_1' \, d\theta. \end{multline*}
Fixing $x$ in the above expression, and integrating out the variables $y_2$, $y_1$, $z_1$, $\theta_1'$, $\theta$ in that order yields the estimate
\begin{equation} 
\sup_x \int \bigl| K_{B_{jk}^{\ast}A_{jk}}(x,y) \bigr| \, dy \leq C 2^{(p+p'+1)j + (l + l'-1)k}. \label{towards-schur} 
\end{equation}    
The same estimate is obtained if the roles of $x$ and $y$ are interchanged in (\ref{towards-schur}). 
Schur's lemma then gives the following estimate for the operator norm:
\begin{eqnarray}\label{old-2}
\quad ||B_{jk}^{\ast} A_{jk}||_{L^2 \rightarrow L^2}  
 & \leq&  \left[\sup_x \int |K_{B_{jk}^{\ast} A_{jk}}(x,y)| dy \right]^{\frac{1}{2}} \left[\sup_y \int |K_{B_{jk}^{\ast} A_{jk}}(x,y)| dx \right]^{\frac{1}{2}}\nonumber \\ 
 &  \leq &C 2^{(l + l'-1)k+(p+ p' + 1)j}.
\end{eqnarray}
Substituting $B = A$, $p=p'$, $l = l'$ yields the conclusion of part (\ref{parta}) of the lemma. 

We now turn to part (\ref{partb}), where we deal with the
``off-diagonal terms'' $B_{j'k'}^{\ast} A_{jk}$ with $|j-j'| \le
3$ and $|k - k'| \geq 3$. Without loss of generality, we may assume
that $k \geq k'+3$. Observing that
\[ (x_1 - y_1) e^{i \tilde{\Phi}} = \bigl[ (z_1 - y_1) + (x_1 - y_1) \bigr] e^{i \tilde{\Phi}} = (d_{\theta_1} + d_{\theta_1'}) e^{i \tilde{\Phi}}, \]
an integration by parts in (\ref{eqn-jkpk}) in the $\theta_1$, $\theta_1'$ variables gives \begin{align*}
\bigl[1 + (2^{k'}(x_1 - y_1))^2 \bigr]^N &\bigl|K_{B_{j'k'}^{\ast} A_{jk}}(x,y) \bigr| \le C_N\times \\ &\left| \int e^{i \tilde{\Phi}(x,y,z; \theta, \theta_1')} u_{jj'kk'}(x,y,z; \theta, \theta')\, dz_1 \, d\theta_1' d \theta \right|.
\end{align*}  
On the other hand, 
\begin{align*}
d_{z_1} \tilde{\Phi} &= \theta_1 - \theta_1' + 3 \theta_2 \bigl((z_1 - x_1)^2 - (z_1 - y_1)^2 \bigr) \\ &= \theta_1 - \theta_1' + 3 \theta_2 \bigl(d^2_{\theta_1'} - d^2_{\theta_1} \bigr) \tilde{\Phi},
\end{align*} 
from which we obtain $L(e^{i \tilde{\Phi}}) = (\theta_1 - \theta_1')
e^{i \tilde{\Phi}}$, where $L := \frac{1}{i} d_{z_1} - 3 \theta_2
d^2_{\theta_1} + 3 \theta_2 d^2_{\theta_1'}$. We observe that $|\theta_1 - \theta_1'| \geq \frac{1}{4} 2^k$, and that the loss from applying
$L^t$ to the amplitude $(\theta_1 - \theta_1')^{-1} u_{jj'kk'}$ is
$\leq C(1 + 2^{j - 2k} + 2^{j - 2k'}) \leq C \max(1, 2^{j - 2k'})$. More precisely, for any $m \geq 1$,  
\[ \min \left(1, 2^{2k'-j} \right) L^t \left[
\frac{u_{jj'kk'}}{(\theta_1 - \theta_1')^m} \right] =
\frac{u_{jj'kk'}}{(\theta_1 - \theta_1')^m}. \] 
Thus an integration by parts argument using the differential operator $(1 + L^2)$ yields
\begin{multline}
\bigl[1 + (2^{k'}(x_1 - y_1))^2 \bigr]^N \bigl|K_{B_{j'k'}^{\ast} A_{jk}}(x,y) \bigr| \le C_N\times \\ \left| \int e^{i \tilde{\Phi}} \left[ 1 + \bigl(\min(1,2^{2k'-j})(\theta_1 - \theta_1') \bigr)^2\right]^{-N} u_{jj'kk'}(x,y,z_1; \theta_2, \theta_1, \theta_1') \, dz_1 \, d\theta_1' \, d\theta \right|. \label{est-comb}
\end{multline}   
We observe that the argument so far does not require any special choice of $\delta$. 

We now consider two cases. First suppose that $k' > \frac{j}{2}$, so that $\min(1, 2^{2k'-j}) = 1$. 
Since $|\theta_1 - \theta_1'| \geq \frac{1}{4}2^k \geq \frac{1}{4}2^{\delta j}$ on the support of $u_{jj'kk'}$, the estimate in (\ref{est-comb}) implies 
\[
\int \bigl| K_{B_{j'k'}^{\ast} A_{jk}}(x,y)\bigr| \, dy + \int \bigl|
K_{B_{j'k'}^{\ast} A_{jk}}(x,y)\bigr| \, dx \leq  C_N 2^{3j - \delta jN} \leq C_N 2^{-jN}
\]
for arbitrarily large $N$. Recalling that $|k-k'| \leq j$, the
inequality above combined with Schur's lemma yields a stronger estimate than the one claimed in (\ref{aok2}).

Next we assume that $2k' \leq j$. Here we will integrate by parts twice more in the integral in (\ref{est-comb}), using the
differential operators $L_r = 1 + \left(\frac{1}{i}
d_{\theta_r}\right)^2$, $r = 1,2$, and keeping in mind that
\[ d_{\theta_1}(e^{i \tilde{\Phi}}) = i (z_1 - y_1)e^{i \tilde{\Phi}},
\quad d_{\theta_2}(e^{i \tilde{\Phi}}) = i(x_2 - y_2 +(z_1 - x_1)^3
-(z_1 - y_1)^3) e^{i \tilde{\Phi}}. \]  We observe that \begin{align*}
d_{\theta_2} (u_{jj'kk'}) &= 2^{-j} u_{jj'kk'}, \text{ while } \\ 
d_{\theta_1} \left[\frac{u_{jj'kk'}}{(1 + (2^{2k'-j}(\theta_1-\theta_1'))^2)^{N}}
 \right] &=  (2^{2k'-j} + 2^{-k})\frac{u_{jj'kk'}}{(1 +
(2^{2k'-j}(\theta_1-\theta_1'))^2)^{N}}, \\ &\leq 2^{2k'-j} \frac{u_{jj'kk'}}{(1 +
(2^{2k'-j}(\theta_1-\theta_1'))^2)^{N}} \end{align*} where the last
inequality follows from the assumptions $k \geq k' \geq
\frac{j'}{3}$ (since $\delta\ge\frac13$) and $|j-j'|\leq 3$. Combining all the arguments above we finally arrive at the estimate    
{\allowdisplaybreaks \begin{multline*}
\bigl[1 + (2^{k'}(x_1 - y_1))^2 \bigr]^N \bigl|K_{B_{j'k'}^{\ast} A_{jk}}(x,y) \bigr| \leq  C_N 2^{(l - \frac{1}{2})k  + (l' - \frac{1}{2})k' + (p+p'+1)j} \times \\ \int \Bigl|\left[ 1 + \bigl(2^{2k'-j}(\theta_1 - \theta_1') \bigr)^2\right]^{-N}  \left[1 + \left( 2^{j - 2k'}|z_1 - y_1|\right)^2 \right]^{-N} \\ \times  \left[1 + \bigl(2^j |x_2 - y_2 + (z_1 - x_1)^3 - (z_1- y_1)^3| \bigr)^2 \right]^{-N} dz_1 \, d\theta_1' \, d\theta \Bigr|. 
\end{multline*}   
Choosing $N$ large enough and observing that $|\theta_1 - \theta_1'|
\geq \frac{1}{4} 2^k$, we obtain 
\begin{equation}
\begin{aligned} \int \bigl|K_{B_{j'k'}^{\ast} A_{jk}} &(x,y) \bigr| \,
dx \\ &\leq \; C 2^{(l - \frac{1}{2})k + (l' - \frac{1}{2})k' +
(p+p'+1)j} 2^{-N(2k'+k-j)} \\ & \quad \times \int_{\theta_2} \int_{\theta_1'} \int_{\theta_1} \left[ 1 +
\bigl(2^{2k'-j}(\theta_1 - \theta_1') \bigr)^2\right]^{-\frac{N}{2}} \\  &
\quad \int_{z_1} \bigl[1 + \bigl( 2^{j - 2k'}|z_1 - y_1|\bigr)^2
\bigr]^{-N} \int_{x_1} \bigl[1 + (2^{k'}(x_1 - y_1))^2 \bigr]^{-N} \\
& \quad \int_{x_2} \frac{dx_2  dx_1  dz_1  d\theta_1
d\theta_1'  d\theta_2}{\left[1 + \bigl(2^j |x_2 - y_2 + (z_1 - x_1)^3 - (z_1-
y_1)^3| \bigr)^2 \right]^{N}}  \\ &\leq  C 2^{(l - \frac{1}{2})k + (l' -
\frac{1}{2}) k' + (p+p'+1)j} 2^{-N(2k' + k - j)},
\end{aligned} \label{IBP-series} \end{equation}}where the last step follows from the fact that the $(x_2, x_1, z_1,
\theta_1, \theta_1', \theta_2)$ integrals when computed in that order
yield $2^{-j}, 2^{-k'}, 2^{2k'-j}, 2^{j - 2k'}, 2^{k'}$ and $2^j$ respectively. By symmetry, the same
estimate also holds for the integral with respect to $dy$ of
$|K_{B_{j'k'}^{\ast} A_{jk}}(x,y)|$, with $x$ fixed. In view of the
last line in (\ref{IBP-series}) and the assumption that $\delta\ge\frac13$, (\ref{aok2}) then follows from Schur's lemma and the inequality \begin{align*} 2k' + k - j =
3k' - j + (k - k') &\geq (3k' - j') + (j - j') + (k-k') \\ & \geq 0 - 2 + (k - k') \end{align*} since $k > k' \geq
\delta j' \geq j'/3$ and $|j-j'|\leq 2$, completing the proof of part (\ref{partb}).   
\end{proof}
\begin{corollary} 
Let $A \in I^{p, l}(\Delta, C_0)$. Then for $j \geq 1$ and $\delta \geq \frac{1}{3}$, the operator $A_j$ satisfies the property
\begin{equation}
 ||A_j||_{L^2 \rightarrow L^2} \leq  \begin{cases}
2^{(p + l)j} &\text{ if } l \geq \frac{1}{2}, \\ 2^{\left[p +
\frac{1}{2} + \delta(l - \frac{1}{2}) \right]j} &\text{ if } l <
\frac{1}{2}. \end{cases} \label{ao-Aj} \end{equation} 
\label{orth-cor}
Further, if $p + l \leq 0$ for $l \geq \frac{1}{2}$, or if $p +
\frac{1}{2} + \delta (l - \frac{1}{2}) \leq 0$ for $l < \frac{1}{2}$,
then $\sum_j A_j = A - A_0 - A_{\infty}$ is a bounded linear map on
$L^2$. 
\end{corollary}
\begin{proof} The proof of (\ref{ao-Aj}) is a direct consequence of
(\ref{aok2}) in Lemma \ref{aok} with $p = p'$ and Cotlar-Knapp-Stein almost orthogonality
lemma \cite{st93}. The second statement uses (\ref{ao-Aj}), Lemma
\ref{aorth1} and almost orthogonality again. The details are left to the reader.
\end{proof}

\section{Normal operators for  seismic imaging}\label{sec-norm}
We now describe the decomposition of an operator $A \in
I^{p,l}(\Delta, C_2)$, where $(\Delta, C_2)$ is as  in \S\S\ref{subsec-marine} and
\S\S\ref{noropseis}. As in \S\ref{sec-gendecomp} let us  
fix $\delta\in[\frac13,\frac12)$, and a dyadic partition of unity $\{
\chi_j : j \geq 0 \}$ on $\mathbb R$ satisfying
(\ref{pou-chi}). Setting \[ \chi_{jk}(\xi, \rho;\sigma) =
\chi_k(\sigma) \chi_j(\langle \xi, \rho \rangle),\;
\chi_{\infty} = 1 - \sum_{j \geq 0} \sum_{k=0}^{j} \chi_{jk}, \text{ where }
\langle \xi, \rho \rangle = (1 + |\xi|^2 + |\rho|^2)^{\frac{1}{2}}, \] 
we arrive at a partition of unity on $\left(\R^{n+1}_{\xi,\rho}\setminus \{0\}\right) \times\R_\sigma$, namely
\[
1\equiv \sum_{j=0}^\infty \sum_{k=0}^j \chi_{jk} + \chi_{\infty} =
\sum_{j=0}^{\infty} \left[ \sum_{k=0}^{[\delta j]} + \sum_{k=[\delta
j]+1}^j \right] \chi_{jk} +\chi_\infty.
\]
Here $\chi_\infty\in S^0_{1,0}$ is supported on $\{ \la\xi,\rho\ra
\leq 2 |\sigma| \}$, while $\{\chi_{jk} : j \geq 0, 0 \leq k \leq j\}$
is a bounded family (with respect to the seminorms) in $S^{0,0}$, with
\begin{equation} \begin{aligned} &\text{supp}(\chi_{jk}) \subseteq\{2^{j-1}\leq \la\xi,\rho\ra\leq 2^{j+1},\, 2^{k-1}
\leq |\sigma| \leq 2^{k+1} \}, \text{ so that } \\
&\text{supp}\Bigl(\sum_{j\geq 0} \sum_{k=0}^{[\delta j]} \chi_{jk} \Bigr) \subseteq
\{ |\sigma|\leq  C\la\xi,\rho\ra^\delta\}, \text{ and } \\ &\text{supp}\Bigl(\sum_{j\geq 0} \sum_{k=[\delta j]+1}^{j} \chi_{jk} \Bigr) \subseteq
\{ C \langle \xi, \rho \rangle^{\delta}  \leq |\sigma|\leq  C\la\xi,\rho\ra\}. \end{aligned} \label{suppjk} \end{equation}  Letting
$a_{jk}(x,y;\xi,\rho;\sigma)=\chi_{jk}(\xi,\rho,\sigma)
a(x,y;\xi,\rho;\sigma)$ and $a_{\infty} = \chi_{\infty}a$ gives rise
to a decomposition of the amplitude $a$, which in turn induces a
decomposition of the operator $A$. More precisely, 
\[ A= A_0 +\sum_{j=0}^{\infty}\sum_{k = [\delta j]+1}^{\infty} A_{jk}+A_{\infty},
\] 
where $A_{jk}$, $A_0$ and $A_{\infty}$ are operators whose Schwartz
kernels and multiphase functions are of the form given in \S\S
\ref{noropseis} (see (\ref{mg-psi}) and (\ref{mg-KA})), but whose amplitudes are
given by $a_{jk}$, $\sum_{j \geq 0} \sum_{0 \leq k \leq [\delta j]}a_{jk}$
and $a_{\infty}$ respectively. Then $A_\infty\in
\Psi^{p+l}$ and so maps $H^s$ boundedly into $H^{s-r}$ for $r\ge
p+l$. It turns out that the components $A_0$ and $A_{jk}$ in this
situation satisfy $L^2$ estimates and almost orthogonality properties
analogous to their $I^{p,l}(\Delta, C_0)$ counterparts with similar but more involved 
proofs. We record these facts below with the appropriate modifications.  

\begin{lemma}[Bounds for $A_0$] \label{boundsA0-C1}
For any $0 < \delta < \frac{1}{2}$, the operator $A_0 \in I^m_{1 -
\delta, \delta}(C_2)$, and hence maps $H^s_{\text{comp}}(\mathbb
R^n)$ boundedly into $H^{s-r_0}_{\text{loc}}(\mathbb R^n)$, where $m$
and $r_0$ are in Lemma \ref{boundsA0}.
\end{lemma} 
\begin{proof}
The argument is identical to the one presented in Lemma \ref{boundsA0}
and involves integrating
out $\sigma$ in the oscillatory representation of $K_{A_0}$. The
details are left to the reader. 
\end{proof} 

\begin{lemma} \label{seismic-aorth1}
Let $A \in I^{p,l}(\Delta, C_2)$, $B \in I^{p', l'}(\Delta, C_2)$,
with decompositions $\{ A_{jk}\}$ and $\{B_{jk} \}$ as above. Then the
conclusions (\ref{aorth1-ineq1}) and (\ref{aorth1-ineq2}) of Lemma
\ref{aorth1} hold for the same set of indices $j$, $j'$, $k$, $k'$ therein.  
\end{lemma}
\begin{proof}
The kernel of $B_{j'k'}^* A_{jk}$ is
\begin{equation}\label{mg-Kjjkk}
\begin{aligned} K_{B_{j'k'}^{\ast}A_{jk}}(x,y)&=\int
e^{i\Phi(x,y;z;\xi,\rho,\sigma;\xt,\rt,\st)} \\ &\qquad \quad \times  
\bar{b}_{j'k'}(z,x;\tilde{\xi}, \tilde{\rho}; \sigma) a_{jk}(z,y;\xi,
\rho; \sigma)  dz d\xi d\rho d\sigma
d\xt d\rt d\st, \end{aligned} 
\end{equation}
where, from (\ref{mg-psi}),
\begin{eqnarray*}
\Phi &= (z-y)\cdot\xi - (z-x)\cdot\xt +\frac{\xi_n\rho}{\xi_1} 
-\frac{\xt_n \rt}{\xt_1} +(z_n-y_n)\sigma - (z_n - x_n)\st \\ &\quad
-\frac14\left( (z_n+y_n)^2\rho - (z_n+x_n)^2\rt\right)
-P(z_n,y,\xi')\frac\rho{\xi_1} +P(z_n,x,\xt')\frac\rt{\xt_1}, 
\end{eqnarray*}
and $\xi' = (\xi_1, \dots, \xi_{n-1})$.
In view of (\ref{xi1-dominant}) and (\ref{suppjk}), the amplitude satisfies the support condition
\[ \text{supp}(\bar{b}_{j'k'} a_{jk})\subseteq \left\{ \begin{aligned} |\langle
\tilde{\xi}, \tilde{\rho} \rangle| \leq 2 |\tilde{\xi}_1|, \quad 2^{j'-1}
\leq |\langle \tilde{\xi}, \tilde{\rho} \rangle|\leq 2^{j'+ 1}, \quad 2^{k'-1} \leq |{\st}|\leq 2^{k'+1} \\
|\langle \xi,\rho \rangle| \leq 2 |\xi_1|, \quad 2^{j-1} \leq |\langle \xi, \rho \rangle| \leq
2^{j+1}, \quad 2^{k-1} \leq |\sigma|\leq 2^{k+1} \end{aligned}
\right\}, \]
and the differentiability estimates
\begin{equation} \label{mg-ajjkk}
\begin{aligned} |\p^\alpha_{\xi,\rho} \p^{\at}_{\xt,\rt}\p^{\beta}_{\sigma} \p^{\bt}_{\st}
\p^\gamma_{x,y,z}(\bar{b}_{j'k'}a_{jk})| &\leq C_{\alpha,
\tilde{\alpha}, \beta, \tilde{\beta}, \gamma}  2^{(p-\frac12)j+(p'-\frac{1}{2})j'
+(l-\frac12)k+(p'-\frac{1}{2})k'} \\
&\hskip1in \times\, 2^{-|\alpha|j-|\at|j'-\beta k -\bt
k'}, \end{aligned} 
\end{equation} 
with $C_{\alpha, \tilde{\alpha}, \beta, \tilde{\beta}, \gamma}$ independent of $j$ and $k$. 

Simply integrating the zeroth order  bounds in (\ref{mg-ajjkk}) yields
the basic estimate
\begin{equation}\label{mg-simple}
\begin{aligned}
|K_{B_{j'k'}^{\ast}A_{jk}}(x,y)| &\leq C 2^{(p-\frac{1}{2})j + (p' -
\frac{1}{2})j' + (l - \frac{1}{2})k + (l' - \frac{1}{2})k'}
2^{(n+1)(j+j')} 2^{(k+k')} \\ 
&\leq C 2^{(p+n + \frac12)(j+j')
+(l+\frac12)(k+k')}.
\end{aligned} 
\end{equation}

Now, if $j'\leq j-3$, then $|\xi_1-\xt_1|\geq c |\xi_1|\geq  
c 2^j$ on supp$(\bar{b}_{j'k'}a_{jk})$. Noting that 
$d_{z_1}\Phi=\xt_1-\xi_1$, we can integrate by parts $N$ 
times in $z_1$ and then estimate as in (\ref{mg-simple}) to obtain
\begin{align*}
|K_{B_{j'k'}^{\ast}A_{jk}}(x,y)| &\leq C_N 2^{-Nj} 2^{(p+n+\frac12)(j+j') 
+(l+\frac12)(k+k')}\\
&\leq C_N  2^{-(N-M)j}
\end{align*}
for some fixed $M = M(p,l,n) \geq 0$, since $k'\le j'<j$ and $k\le
j$. By choosing $N$ sufficiently large compared to $M$, integrating 
the kernel in $x$ and taking the supremum in $y$, or vice versa, we
get the same upper bound, which establishes the desired $L^2$ norm by
Schur's lemma. The case of $j' \geq j+3$ is identical.
\end{proof} 
\begin{lemma}[Almost orthogonality in $k$]
Let $A$, $B$, $\{ A_{jk}\}$ and $\{B_{jk} \}$ be as in Lemma
\ref{seismic-aorth1}. If the indices $j$, $j'$, $k$, $k'$ satisfy the same
hypotheses as in Lemma \ref{aok}, then the conclusions (\ref{AjkL2}) and (\ref{aok2}) hold.  
\end{lemma}
\begin{proof}
The proof of (\ref{AjkL2}) and (\ref{aok2}) in this setting requires an initial
preparation of the Schwartz kernel of $B_{j'k'}^{\ast}A_{jk}$; we
simplify (\ref{mg-Kjjkk}) using four applications of stationary
phase. The first three applications are with respect to the pairs
$(\rho,\xi_n)$, $(\rt,\xt_n)$ and $(z',\tilde{\xi}')$. Of these the first
method of stationary phase may be justified using Lemma \ref{SP-lemma}, by setting $w = (\rho, \xi_n)$,
so that in the notation of that lemma \[ Q = \left[ \begin{matrix} 0 & \frac{1}{\xi_1} \\
\frac{1}{\xi_1} & 0  \end{matrix} \right], \quad \lambda^{-1}R = -2\xi_1
\frac{\partial^2}{\partial \xi_n \partial \rho}, \quad \text{ and } \quad
(\lambda^{-1}R)^r(\bar{b}_{j'k'}a_{jk}) = O(2^{-jr}), \] 
giving rise to a valid asymptotic expansion whose first term involves
a factor $\det(Q)^{-\frac{1}{2}} = |\xi_1| \sim 2^j$. The other stationary
phases are handled similarly, their justification being left to the
reader. By the same argument as before, the second stationary phase contributes a factor of
$2^{j'}$ (from the determinant of the Hessian) to the
leading order term of the asymptotic expansion, while the
corresponding contribution from the third is only a constant. The end result is
\begin{equation}\label{K3}
K_{B_{j'k'}^{\ast}A_{jk}}(x,y) \sim\int
e^{i\Phi_3} c_3 dz_n  d\xi' d\sigma d\st,
\end{equation}
where $\Phi_3$ and $c_3$ are both functions of
$(x,y;z_n;\xi',\sigma,\st)$, with 
\begin{eqnarray*}\label{Phi3}
\quad\Phi_3= (x'-y')\cdot\xi' & + & (z_n-y_n)\sigma - (z_n-x_n)\st\\
& - & \frac14\left( 
(z_n-x_n)(z_n+x_n)^2-(z_n-y_n)(z_n+y_n)^2\right)\xi_1\nonumber\\
&+& (z_n-y_n)P(z_n,y,\xi')-(z_n-x_n)P(z_n,x,\xi')\nonumber, \text{ and
} 
\end{eqnarray*}
\begin{align*}
c_3 &= \xi_1^2\times \bar{b}_{j'k'}\bigl(x', z_n, x; \xi', \frac{1}{4}(z_n +
x_n)^2 \xi_1 + P(z_n, x, \xi'); (x_n -z_n) \xi_1, \tilde{\sigma}
\bigr) \\ &\qquad \times a_{jk}\bigl(x', z_n, y; \xi', \frac{1}{4}(z_n
+ y_n)^2 \xi_1 + P(z_n, y, \xi'); (y_n - z_n)\xi_1, \sigma \bigr)
\end{align*} 
The amplitude $c_3$ satisfies the same size and differentiability estimates as
$2^{j+j'}\bar{b}_{j'k'}a_{jk}$ with the eliminated variables
absent (see (\ref{mg-ajjkk})). 

To estimate the norm of each $\ajk$, we
set $B = A$, $j=j'$, $k=k'$, $p=p'$, $l=l'$ and apply one final method of stationary phase in
$(z_n,\st)$ (once again justified by the existence of the large
parameter $|\tilde{\sigma}| \geq 2^{k-1} \gg 1$), obtaining 
\begin{equation}\label{K4}
K_{A_{jk}^{\ast}A_{jk}}(x,y)\sim \int e^{i\Phi_4} c_4 d\xi' d\sigma.
\end{equation}
Here $\Phi_4$ and $c_4$ are functions of $(x,y;\xi', \sigma)$, with 
\begin{eqnarray}\label{Phi4}
\Phi_4 &=& (x'-y')\cdot\xi' +(x_n-y_n)\sigma 
+\frac14(x_n-y_n)(x_n+y_n)^2\xi_1\\
& & + (x_n-y_n) P(x_n,y,\xi')\nonumber,
\end{eqnarray}
and $c_4$ satisfying the same estimates as $c_3$, i.e.,
\begin{equation}
|\p^\alpha_{\xi'} \p^\beta_\sigma \p^\gamma_{x,y} c_4|\leq C_{\alpha
\beta \gamma} 2^{(2p+1)j+(2l-1)k
-|\alpha|j-|\beta|k}. \label{c4est} 
\end{equation}
Simple integration in $\xi', \sigma$ in the range $|\xi'| \leq
2^{j+1}, |\sigma| \leq 2^{k+1}$ yields the basic estimate 
\[\left|K_{A_{jk}^{\ast}A_{jk}}(x,y)\right|\leq C 2^{(2p+n)j+2lk},\] but
one may substantially improve upon this by integrating by parts, using
the differential operators 
\begin{eqnarray}\label{3derivs}
\quad d_{\xi_1}\Phi_4 &=& x_1-y_1+\frac14(x_n-y_n)(x_n+y_n)^2 +(x_n-y_n)
P_{\xi_1}(x_n,y,\xi')\nonumber, \\ \quad d_{\xi''}\Phi_4 &=&
x''-y''+(x_n-y_n)P_{\xi''}(x_n,y,\xi'), \\ \quad d_\sigma\Phi_4 &=& x_n-y_n,\nonumber
\end{eqnarray}
where $x''=(x_2,\dots,x_{n-1})$. Since $\partial_{\xi'}^{\beta}P(x_n,
y, \xi')$ is homogeneous of degree $1 - |\beta|$ in $\xi'$, and
$|\xi'| \geq 2^{j-1}$, we find that 
\begin{equation} \label{ibp-justify}
\begin{aligned}
\frac{\partial}{\partial \xi_i} \left[\frac{\tilde{c}_4}{(d_{\xi_i}\Phi_4)^m}
\right] &= 2^{-j} \frac{\tilde{c}_4}{(d_{\xi_i}\Phi_4)^m}, \quad 1 \leq i
\leq n-1, \, m \geq 1, \quad \text{ while } \\ 
\frac{\partial}{\partial \sigma}(\tilde{c}_4) &= 2^{-k} \tilde{c}_4,
\end{aligned}
\end{equation} 
where $\tilde{c}_4$ denotes a function satisfying the same estimates
(\ref{c4est}) as $c_4$, possibly with different implicit constants, and
whose exact functional form may vary from one occurrence to the next.  
Integrating by parts $N$ times with respect to each of the variables
$(\xi_1, \dots, \xi_{n-1}, \sigma)$ in (\ref{K4}), and combining
(\ref{3derivs}) and (\ref{ibp-justify}) yields
\begin{equation}\nonumber
\left|K_{A_{jk}^{\ast}A_{jk}}(x,y)\right| \leq C_N \int 
\frac{2^{(2p+1)j+(2l-1)k}\quad d\xi'd\sigma}
{(1+2^j|d_{\xi'}\Phi_4|)^N 
(1+2^k|d_\sigma\Phi_4|)^N},
\end{equation}
from which one obtains
\begin{align*}
\int|K_{A_{jk}^{\ast}A_{jk}}(x,y)|dx &\leq C_N 
2^{(2p+1)j+(2l-1)k}\iint \frac{d\xi' d\sigma dx}
{(1+2^j|d_{\xi'}\Phi_4|)^N 
(1+2^k|d_\sigma\Phi_4|)^N}\\
&\leq C_N 2^{(2p+n)j+(2l-1)k} 2^{-(n-1)j-k}\int_{\begin{subarray}{c}
|\xi'| \leq 2^{j+1} \\ |\sigma| \leq 2^{k+1} \end{subarray}} d\xi' d\sigma\\
&\leq C 2^{(2p+1)j+(2l-1)k}.
\end{align*}
In the second step above and in view of (\ref{3derivs}), one has to
integrate first in $x'$ and then in $x_n$, while the last step uses
the size of the $(\xi', \sigma)$-support of $a_{jk}$. 
The same estimate holds for $\int|K_{A_{jk}^{\ast}A_{jk}}(x,y)| dy$,
completing the proof of (\ref{AjkL2}) for $I^{p,l}(\Delta, C_2)$.

In order to prove the analogue of (\ref{aok2}), it suffices to show that for any $N \geq 1$, $|j-j'|\leq 2$, $|k-k'|\geq 3$, 
\begin{multline}
\sup_{x}\int \left| K_{B_{j'k'}^{\ast}A_{jk}}(x,y) \right| \, dy + \sup_{y} \int \left| K_{B_{j'k'}^{\ast}A_{jk}}(x,y) \right|\, dx \\ \leq C_N 2^{(p+p'+1)j + (l - \frac{1}{2})k + (l'- \frac{1}{2})k'} 2^{-|k - k'|N}.  \label{aok-C1} 
\end{multline} 
For this we return to the
representation (\ref{K3}). As in the proof of Lemma \ref{aok}, we will
subject the integral in (\ref{K3}) to a large of number of integration
by parts using several differential operators. However due to the complicated structure of the multiphase function (\ref{Phi3}) (compared to (\ref{eqn-Phit})), the differential operators and hence the resulting integrations by parts are more involved. We now proceed to describe each of these steps systematically. Throughout this discussion, as in the proof of Lemma \ref{aok}, we will assume that $k \geq k' + 3$, and denote by $u_{jj'kk'}$ any function (with possibly different functional forms) satisfying the same size and differentiability estimates as $c_3$. The value of the large constant $N$ may also vary from one occurrence to the next.    

{\bf{Step 1. }} Since
\begin{equation} d_\sigma\Phi_3=z_n-y_n,\quad d_{\tilde\sigma}\Phi_3=x_n-z_n,\quad
(d_\sigma+d_{\st})\Phi_3=x_n-y_n, \label{ibp-s-st} \end{equation} 
with \[d_{\sigma}(u_{jj'kk'}) = 2^{-k} u_{jj'kk'} \quad \text{ and } \quad d_{\tilde{\sigma}} (u_{jj'kk'}) = 2^{-k'}u_{jj'kk'}, \] integrating by parts in (\ref{K3}) a large number of times using $1 + \left[\frac{1}{i}(d_{\sigma} + d_{\tilde{\sigma}})\right]^2$ and $1 + (\frac{1}{i}d_{\tilde{\sigma}})^2$ gives  
\begin{equation} \label{IBP-step1}
\begin{aligned} \left[ 1 + \bigl(2^{k'}(x_n - y_n)\bigr)^2 \right]^{N} &\left| K_{B_{j'k'}^{\ast}A_{jk}}(x,y)\right| \\  \sim &\left| \int e^{i \Phi_3} \frac{u_{jj'kk'}(x,y,z_n; \sigma, \tilde{\sigma})}{\left[1 + (2^{k'}(x_n - z_n))^2 \right]^N}\, dz_n \, d\xi' \, d\sigma \, d\tilde{\sigma} \right|. \end{aligned} 
\end{equation} 
The above expression permits a localization in the spatial variables $x_n, y_n, z_n$, which will be useful in the sequel. Fixing $\chi \in C_0^{\infty}([-2,2])$, with $\chi \equiv 1$ on $[-1,1]$, and introducing the partition of unity 
\begin{align*} 1 &\equiv \chi(2^{k'-k' \epsilon}(x_n - y_n)) + (1 - \chi)(2^{k'-k'\epsilon}(x_n - y_n)) \\ &\equiv \chi(2^{k'-k' \epsilon}(x_n - y_n))\left[ \chi(2^{k'-k'\epsilon}(x_n - z_n)) + (1 - \chi)(2^{k'-k'\epsilon}(x_n - z_n)) \right] \\ &\hskip 3in + (1 - \chi)(2^{k'-k'\epsilon}(x_n - y_n)),  \end{align*}  
we obtain the decomposition \[ K_{B_{j'k'}^{\ast}A_{jk}} = \mathcal K_1 + \mathcal K_2 + \mathcal K_3, \] 
where each $\mathcal K_i$ is an oscillatory integral with multi-phase $\Phi_3$, and amplitude of the form $(1 + (2^{k'}(x_n - z_n))^2)^{-N} (1 + (2^{k'}(x_n - y_n))^2)^{-N} \Theta_i$, where 
\begin{align*}
\Theta_1 &= \chi(2^{k'-k' \epsilon}(x_n - y_n))\chi(2^{k'-k'\epsilon}(x_n - z_n)) u_{jj'kk'}, \\
\Theta_2 &= \chi(2^{k'-k' \epsilon}(x_n - y_n))(1 - \chi)(2^{k'-k'\epsilon}(x_n - z_n)) u_{jj'kk'}  \\ 
\Theta_3 &= (1 - \chi)(2^{k' - k'\epsilon}(x_n - y_n)) u_{jj'kk'}. 
\end{align*} 
It follows from (\ref{IBP-step1}) that for any $\epsilon > 0$
\[ \left| \mathcal K_i(x,y) \right| \leq C_N 2^{-jN} \quad \text{ for } i = 2,3 \text{ and all } N \geq 1, \]
which is a stronger statement than the one required in (\ref{aok-C1}). We therefore restrict attention only to $\mathcal K_1$ in the sequel, in which $x_n, y_n$ and $z_n$ are further restricted to satisfy 
\begin{equation} |x_n - z_n| \leq 2^{1-k' + k' \epsilon}, \; |x_n - y_n| \leq 2^{1-k' + k' \epsilon}, \text{ hence } |z_n - y_n| \leq C2^{-k' + k'\epsilon}.  \label{xnynzn-supp} \end{equation}     

{\bf{Step 2. }} Our next integration by parts (in $\mathcal K_1$) will involve $z_n$, $\sigma$, $\tilde{\sigma}$ and will exploit the disparity in the sizes of $\sigma$ and $\tilde{\sigma}$. For this, we note that   
\begin{align*}
d_{z_n}\Phi_3 = \sigma - \tilde{\sigma} &+ \frac{\xi_1}{4} \left[2(x_n - z_n)(x_n - y_n) - (x_n - y_n)^2\right] + P(z_n, y, \xi') - P(z_n, x, \xi') \\ &+ (z_n - y_n) P_{z_n}(z_n, y, \xi') - (z_n - x_n) P_{z_n}(z_n, x, \xi') \\ = \sigma - \tilde{\sigma} &+ \frac{\xi_1}{4} \left[2(x_n - z_n)(x_n - y_n) - (x_n - y_n)^2\right] \\ &+ \frac{3}{2}(z_n - y_n)^2 P_{z_nz_n} (y_n, y, \xi') - \frac{3}{2} (z_n - x_n)^2 P_{z_n z_n} (x_n, x, \xi') \\ &+ (z_n-y_n)^3 \mathcal Q(z_n, y, \xi') - (z_n - x_n)^3 \mathcal Q(z_n, x, \xi'),
\end{align*} 
where the last step follows by expanding $P(z_n, y, \xi')$ and $P_{z_n}(z_n, y, \xi')$ (resp.\ $P(z_n, x, \xi')$ and $P_{z_n}(z_n, x, \xi')$) in a Taylor series in the $z_n$ variable about $z_n = y_n$ (resp.\ $z_n = x_n$) using (\ref{mg-P}). Here $\mathcal Q$ is a smooth function that is homogeneous of degree one in $\xi'$. In view of (\ref{ibp-s-st}), we find that  
\begin{align*}
L(\Phi_3) &= \sigma - \tilde{\sigma} + (z_n - y_n)^3 \mathcal Q(z_n, y, \xi') - (z_n - x_n)^3 \mathcal Q(z_n, x, \xi'), \text{ where } \\ 
L &= d_{z_n} - \frac{\xi_1}{4} \left[ 2 d_{\tilde{\sigma}}(d_{\sigma} + d_{\tilde{\sigma}}) - (d_{\sigma} + d_{\tilde{\sigma}})^2\right] \\ &\hskip1.5in - \frac{3}{2} P_{z_nz_n}(y_n, y, \xi') d_{\sigma}^2 + \frac{3}{2} P_{z_nz_n}(x_n, x, \xi') d^2_{\tilde{\sigma}} \\ &= L_1 + L_2, \text\quad { with } \quad L_1 := d_{z_n}, \text{ and } \\ L_2 &= \mathcal P(y_n, y, \xi') d_{\sigma}^2 - \mathcal P(x_n, x, \xi') d_{\tilde{\sigma}}^2.  
\end{align*} 
Here $\mathcal P(x_n, x, \xi') = \frac{\xi_1}{4} - \frac{3}{2}P_{z_nz_n}(x_n, x, \xi')$ is a smooth function in all its arguments and homogeneous of degree one in $\xi'$. The support properties in (\ref{xnynzn-supp}) and the homogeneity of $\mathcal Q$ imply that 
\begin{equation} \label{ibp-L}
\begin{aligned} 
|L(\Phi_3)| &\geq |\sigma| - |\tilde{\sigma}| - C2^j|z_n - y_n|^3 - C2^j|z_n - x_n|^3 \\ &\geq 2^{k-1} - 2^{k'+1} - C2^{j-3k+3k \epsilon} - C2^{j-3k'+3k' \epsilon} \geq c 2^k, 
\end{aligned} 
\end{equation}   
where the last inequality follows from $k \geq k' + 3$ and $k + 3k'(1 - \epsilon) - j \geq j'(1 - \epsilon) - \frac{2j}{3} \geq j(\frac{1}{3} - \epsilon) - 2 \gg 1$, choosing $\epsilon < \frac{1}{3}$.    

We will integrate by parts a large number of times in $\mathcal K_1$ using $L$, which is justified in light of (\ref{ibp-L}). In order to describe the action of the differential operator on the amplitude, let us denote by $\tilde{\Theta}_1$ any (generic) function satisfying the same size and differentiability estimates as $\Theta_1$, and observe that $L_2$ is independent of $z_n$, so that for any $m \geq 1$, \begin{equation} \label{L2action} \begin{aligned} L_2^t\left[ \frac{\tilde{\Theta}_1}{L(\Phi_3)^m}\right] &=  \sum_{r = 0}^{2} \frac{1}{L(\Phi_3)^{m+r}}\left[\mathcal P(y_n, y, \xi') d_{\sigma}^{2-r} - \mathcal P(x_n, x, \xi')  d_{\tilde{\sigma}}^{2-r} \right](\tilde{\Theta}_1) \\ &= 2^{j} \sum_{r=0}^{2} \frac{2^{-k'(2-r)}\tilde{\Theta}_1}{(L(\Phi_3))^{m+r}}, \quad \text{ since } k \geq k'+3 \text{ and } \mathcal P \tilde{\Theta}_1 = 2^j \tilde{\Theta}_1. \end{aligned} \end{equation} 
On the other hand, \begin{align*}L_1^t &\left[ \frac{\tilde{\Theta}_1}{(L(\Phi_3))^m \bigl(1 + (2^{k'}(x_n - z_n))^2 \bigr)^N} \right] =  \frac{d_{z_n} \tilde{\Theta}_1}{(L(\Phi_3))^m \left[ 1 + (2^{k'}(x_n - z_n))^2\right]^N} \\ & \hskip1in + \frac{\tilde{\Theta}_1}{(L(\Phi_3))^{m+1}} \left[ \frac{(z_n - y_n)^2 \tilde{\mathcal Q}(z_n, y, \xi') - (z_n - x_n)^2 \tilde{\mathcal Q}(z_n, x, \xi')}{\left[ 1 + (2^{k'}(x_n - z_n))^2\right]^N}\right] \\ &\hskip2in +\frac{\tilde{\Theta}_1}{(L(\Phi_3))^m} \frac{2^{2k'}(x_n - z_n)}{\left[ 1 + (2^{k'}(x_n - z_n))^2\right]^{N+1}} \\ &\hskip1.9in = \frac{2^{k'-k'\epsilon}\tilde{\Theta}_1}{(L(\Phi_3))^m \left[ 1+(2^{k'}(x_n - z_n))^2\right]^N} \\ &\hskip1in+ \frac{1}{(L(\Phi_3))^{m+1}} \left[ \frac{2^j(z_n - y_n)^2 \tilde{\Theta}_1 - 2^j(z_n - x_n)^2 \tilde{\Theta}_1}{\left[ 1 + (2^{k'}(x_n - z_n))^2\right]^N}\right] \\ &\hskip2in +\frac{\tilde{\Theta}_1}{(L(\Phi_3))^m} \frac{2^{2k'}(x_n - z_n)}{\left[ 1 + (2^{k'}(x_n - z_n))^2\right]^{N}} \end{align*}  
where $\tilde{\mathcal Q}$ shares the same smoothness and homogeneity properties as $\mathcal Q$, and hence $\tilde{\Theta}_1 \tilde{\mathcal Q} = 2^j \tilde{\Theta}_1$. In view of (\ref{ibp-s-st}) and the factors involving $(z_n - y_n), (x_n - z_n)$ in the expression above, one can follow up an application of $L_1^t$ by another integration by parts using $d_{\sigma}$ and $d_{\tilde{\sigma}}$, obtaining 
\begin{equation} \label{L1action} \begin{aligned} 
\int &e^{i \Phi_3} L_1^t \left[ \frac{\tilde{\Theta}_1}{(L(\Phi_3))^m \bigl(1 + (2^{k'}(x_n - z_n))^2 \bigr)^N} \right] \, dz_n \, d\xi' \, d\sigma \, d\tilde{\sigma} \\  &= \int \frac{e^{i \Phi_3}}{\Bigl[ 1 + (2^{k'}(x_n - z_n))^2 \Bigr]^N}\Biggl[\frac{2^{k'-k'\epsilon}}{(L(\Phi_3))^m} + 2^j (d_{\sigma}^{2} - d^2_{\tilde{\sigma}})\left(\frac{\tilde{\Theta}_1}{(L(\Phi_3))^{m+1}} \right) \\ &\hskip2in + 2^{2k'} d_{\tilde{\sigma}}\left(\frac{\tilde{\Theta}_1}{(L(\Phi_3)^m)}\right) \Biggr] \, dz_n \, d\xi' \, d\sigma \, d\tilde{\sigma} \\ &= \int \frac{e^{i \Phi_3}}{\Bigl[ 1 + (2^{k'}(x_n - z_n))^2 \Bigr]^N}\Biggl[ \frac{2^{k'} \tilde{\Theta}_1}{(L(\Phi_3))^m} + \frac{2^{2k'} \tilde{\Theta}_1}{(L(\Phi_3))^{m+1}} \\ &\hskip2in + 2^j \sum_{r=0}^{2} \frac{2^{-k'(2-r)} \tilde{\Theta}_1}{(L(\Phi_3))^{m+1+r}} \Biggr] \, dz_n \, d\xi' \, d\sigma d\tilde{\sigma}  
\end{aligned} \end{equation}   
Combining (\ref{L2action}) and (\ref{L1action}), we observe that a $J$-fold application of integration by parts in $\mathcal K_1$ using $L$ results in a finite sum of oscillatory integrals with multiphase $\Phi_3$, whose amplitudes are all of a similar form, namely
\begin{equation} \label{amplitude-form}
\Xi_1  = \frac{2^{\beta_J(j,k', \mu)}\tilde{\Theta}_1}{(L(\Phi_3))^{\mu} \left[ 1 + (2^{k'}(x_n - z_n))^2\right]^{N} \left[1 + (2^{k'}(x_n - y_n))^2 \right]^N},
\end{equation}  
with $\mu \geq J$. It is easy to check using the restrictions $3k' \leq j$ and $k' \leq k$ that the exponents $\beta_1(j,k',\mu)$ satisfy the inequality $\beta_1(j,k',\mu) \leq (\mu-1) k + k'$. An induction in $J$, which we ask the reader to verify, then shows that in general 
\begin{equation} \beta_J(j,k',\mu) \leq (\mu - J) k + k'J, \quad J\geq 1. \label{kappa-ineq} \end{equation} We work with a general term in this sum, which by a slight abuse of notation we continue to denote by $\mathcal K_1$:   
\begin{equation} \mathcal K_1(x,y) \sim \int e^{i \Phi_3} \Xi_1(x,y,z_n; \xi', \sigma, \tilde{\sigma})\, dz_n \, d\xi' \, d\sigma \, d\tilde{\sigma}, \label{K1-revised} \end{equation} 
where $\Xi_1$ is as in (\ref{amplitude-form}).

{\bf{Step 3. }} Our last integration by parts will be in the variables $\xi', \sigma, \tilde{\sigma}$. Specifically, we note that  
\begin{align*}
d_{\xi'}\Phi_3 & -\left[x'-y'+ \frac{1}{4}\left[(z_n-y_n)(z_n+y_n)^2-(z_n-x_n)(z_n+x_n)^2 \right] \vec{e}_1 \right] \\
& \quad = (z_n-y_n)P_{\xi'}(z_n,y,\xi')-(z_n-x_n) P_{\xi'}(z_n,x,\xi') \\ &\quad = \mathcal R(z_n, y, \xi') (z_n - y_n)^3 - (z_n - x_n)^3 \mathcal R(z_n, x, \xi'),
\end{align*}
where the last step follows from Taylor expansion of $P_{\xi'}(z_n, y, \xi')$ (resp.\ $P_{\xi'}(z_n, x, \xi')$) in the $z_n$ variable about the point $z_n = y_n$ (resp.\ $z_n = x_n$) using (\ref{mg-P}). Here $\mathcal R = (\mathcal R_1, \cdots, \mathcal R_{n-1})$ is a smooth vector-valued function that is homogeneous of degree zero in $\xi'$. In view of (\ref{ibp-s-st}), this implies that the vector differential operator $\mathbf P = (P_1, \cdots, P_{n-1}) = d_{\xi'} - \mathcal R(z_n, y, \xi') d_{\sigma}^3 - \mathcal R(z_n, x, \xi') d_{\tilde{\sigma}}^3$ satisfies 
\begin{equation}\label{P} \mathbf P(\Phi_3) = x'-y'+ \frac{1}{4}\left[(z_n-y_n)(z_n+y_n)^2-(z_n-x_n)(z_n+x_n)^2 \right] \vec{e}_1.  \end{equation}
In order to describe the result of a large number of integration by parts in $\mathcal K_1$ using $\mathbf P$, we observe that $\mathbf P(\Phi_3)$ is independent of $(\xi',\sigma,\tilde{\sigma})$, so (in view of (\ref{K1-revised}) and (\ref{amplitude-form})) one only needs to understand the effect of $\mathbf P^t$ on $\tilde{\Theta}_1 (L(\Phi_3))^{-\mu}$. For every $1 \leq i \leq n-1$,   
\begin{align*}
P_i^t \left[\frac{\tilde{\Theta}_1}{(L(\Phi_3))^{\mu}} \right] = & \frac{d_{\xi_i} \tilde{\Theta}_1}{(L(\Phi_3))^{\mu}} \\ &+ \frac{\tilde{\Theta}_1}{(L(\Phi_3))^{\mu+1}} \left[(z_n - y_n)^3 \mathcal Q_{\xi_i}(z_n, y, \xi') - (z_n - x_n)^3 \mathcal Q_{\xi_i}(z_n, x, \xi') \right] \\
& - \left[\mathcal R(z_n, x, \xi') d_{\tilde{\sigma}}^3 + \mathcal R(z_n, x, \xi') d_{\tilde{\sigma}}^3 \right] \left( \frac{\tilde{\Theta}_1}{(L(\Phi_3))^{\mu}}\right)
\end{align*}
As in step 2, we exploit (\ref{ibp-s-st})  and use another integration by parts with respect to $\sigma$, $\tilde{\sigma}$ to replace the factors involving $(z_n - y_n)$ and $(x_n - z_n)$ in the expression above, thus obtaining
\begin{equation} \label{Pi-calc}
\begin{aligned}
\int &e^{i \Phi_3} P_i^t\left[ \frac{\tilde{\Theta}_1}{\left( L(\Phi_3)\right)^{\mu}} \right] d\xi' \, d\sigma \, d\tilde{\sigma} = 
\int d\xi' \, d\sigma \, d\tilde{\sigma} e^{i \Phi_3} \Biggl[ \frac{2^{-j} \tilde{\Theta}_1}{\left( L(\Phi_3)\right)^{\mu}} \\ &\hskip1in + \left[\mathcal Q_{\xi_i}(z_n, y, \xi') d_{\sigma}^3 + \mathcal Q_{\xi_i}(z_n, x, \xi') d^3_{\tilde{\sigma}} \right] \left( \frac{\tilde{\Theta}_1}{\left( L(\Phi_3)\right)^{\mu+1}}\right) \\ & \hskip1in - \left[\mathcal R_i(z_n, y, \xi') d_{\sigma}^3 + \mathcal R_i(z_n, x, \xi')d^3_{\tilde{\sigma}} \right] \left(\frac{\tilde{\Theta}_1}{\left(L(\Phi_3) \right)^{\mu}} \right) \Biggr] \\ &= \int e^{i \Phi_3} \Biggl[ \frac{2^{-j} \tilde{\Theta}_1}{\left(L(\Phi_3)\right)^{\mu}} + \sum_{r=0}^{3} \frac{2^{-k'(3-r)} \tilde{\Theta}_1}{\left( L(\Phi_3 )\right)^{\mu + 1 + r}} + \sum_{r=0}^{3} \frac{2^{-k'(3-r)} \tilde{\Theta}_1}{\left( L(\Phi_3 )\right)^{\mu + r}} \Biggr]\, d\xi' \, d\sigma \, d\tilde{\sigma}, 
 \end{aligned} 
\end{equation} 
since $\mathcal Q_{\xi_i} \tilde{\Theta}_1 = \tilde{\Theta}_1$ and $\mathcal R_i \tilde{\Theta}_1 = \tilde{\Theta}_1$. Recalling (\ref{ibp-L}) and observing that 
\[ k'(3-r) + k(1 + r) \geq k'(3-r) + kr \geq 3k' \geq j \] 
we conclude from (\ref{Pi-calc}) after an easy induction that an $m$-fold integration by parts using $P_i$ yields 
\begin{equation} \int e^{i \Phi_3} \frac{\tilde{\Theta}_1}{L(\Phi_3)^{\mu}} \, d\xi' \, d\sigma \, d\tilde{\sigma} = \frac{1}{\left( P_i \Phi_3 \right)^m}\sum_{r \geq 0} 2^{\gamma_{m}(j,k',r)} \int e^{i \Phi_3} \frac{\tilde{\Theta}_1 \, d\xi' d\sigma d\tilde{\sigma}}{\left( L(\Phi_3)\right)^{\mu + r}}\label{ibp-final}  \end{equation}
where \begin{equation} \label{nu-ineq} \gamma_m(j,k',r) \leq -jm. \end{equation} This concludes the steps that involve various integration by parts.

{\bf{Step 4.}} It remains to combine the results of the previous steps to prove (\ref{aok-C1}). Combining (\ref{K1-revised}) with (\ref{ibp-final}), choosing $m=M \gg 1$ and continuing to call a generic term in the resulting sum by $\mathcal K_1$, we find that 
\begin{multline*} \mathcal K_1(x,y) \le C_N \frac{2^{\beta_J(j,k',\mu)}}{\left[1 + (2^{k'}(x_n - y_n))^2\right]^N} \\ \times \int \frac{e^{i \Phi_3} \tilde{\Theta}_1 \, dz_n d\xi' d\sigma d\tilde{\sigma}}{\left( L(\Phi_3)\right)^{\mu + r}  \left[ 1 + \left(2^{k'}(z_n - x_n)\right)^2\right]^N  \prod_{i=1}^{n-1}\left[1 + 2^{-\gamma_{M}(j,k',r)} (P_{i}\Phi_3)^M \right]}. \end{multline*}
To prove (\ref{aok-C1}) we compute 
\begin{align*} &\sup_x \int |\mathcal K_1(x,y)| \, dy \leq C2^{(p + p' + 1)j + (l - \frac{1}{2})k + (l' - \frac{1}{2})k'} 2^{\beta_J(j,k',\mu) - k\mu} \times \\ & \sup_x \int \frac{dy' \, dz_n \, dy_n \, d\xi' d\sigma \, d\tilde{\sigma}}{\left[1 + \left( 2^{k'}(x_n - y_n)\right)^2 + \left( 2^{k'}(z_n - y_n)\right)^{2} \right]^N \prod_{i=1}^{n-1}\left[1 + 2^{-\gamma_{M}(j,k'r)} (P_{i}\Phi_3)^M \right]} \\ &\leq C2^{(p + p' + 1)j + (l - \frac{1}{2})k + (l' - \frac{1}{2})k'} 2^{\beta_J(j,k',\mu) - k\mu} 2^{(n-1)\frac{\gamma_M(j,k',r)}{M}} 2^{-k'} 2^{-k'} 2^{j(n-1)} 2^k 2^{k'} \\ &\leq C2^{(p + p' + 1)j + (l - \frac{1}{2})k + (l' - \frac{1}{2})k'}2^{\beta_J(j,k',\mu) - k\mu + k - k'} \\ & \leq C 2^{(p + p' + 1)j + (l - \frac{1}{2})k + (l' - \frac{1}{2})k'} 2^{-(k-k')(J-1)}, \end{align*}  
where we have used the size estimate for $\tilde{\Theta}_1$, $r\geq 0$ and (\ref{ibp-L}) at the first step; the expressions (\ref{P}) and the support sizes of $a_{jk}$, $b_{j'k'}$ at the second step to integrate in the order specified; (\ref{nu-ineq}) at the third step; and (\ref{kappa-ineq}) at the final step. Since the roles of $x$ and $y$ may be interchanged to get the same bound, the desired conclusion is established.    
\end{proof}

\begin{corollary} \label{orth-cor-C1}
For $\delta \geq \frac{1}{3}$, the statement of Cor.\ \ref{orth-cor} holds for $A \in I^{p,l}(\Delta, C_2)$. 
\end{corollary}
\begin{proof}
The proof is identical to that of Cor.\ \ref{orth-cor}. 
\end{proof}

\section{Proof of Theorem \ref{thm-main}}\label{sec-cubic}
Now in a position to prove Thm.\ \ref{thm-main},  we first treat  the classes of operators $I^{p,l}(\Delta, C_j)$ for $j=0,2$.  
Since there are standard elliptic  $\Psi$DOs that map any $L^2$-Sobolev space isomorphically  onto any other, and
\[ \Psi^m \circ I^{p, l}(\Delta, C_j) \circ \Psi^{m'} \subseteq I^{p + m + m', l}(\Delta, C_j) \quad \text{for any }\, m,m'\in\R, \] 
it suffices to assume that $s = r = 0$ in Thm. \ref{thm-main} and  show that $A \in I^{p,l}(\Delta, C_j)$ is  bounded  from $L^2 \rightarrow L^2$. These assumptions will be used without further reference throughout this section as will the decompositions of the operator $A$ introduced in $\S\ref{sec-gendecomp}$ and $\S\ref{sec-norm}$.       
The proof consists of combining the estimates for $A_0$, $A_{\infty}$ and $A - A_0 - A_{\infty}$. Of these $A_{\infty} \in \Psi^{p+l}$, and hence is a bounded map from $H^s \rightarrow H^{s-r'}$, for all $s \in \mathbb R$ and $r' \geq r_{\infty} = p+l$. For $I^{p,l}(\Delta, C_0)$, the bounds for $A_0$ and $A - A_0 - A_{\infty}$ are given in Lemma \ref{boundsA0} and Cor.\ \ref{orth-cor} respectively; for $I^{p,l}(\Delta, C_2)$ they are in Lemma \ref{boundsA0-C1} and Cor.\ \ref{orth-cor-C1}.  
It therefore suffices to show that \[r =  \min_{\delta \in [\frac{1}{3}, \frac{1}{2})} \max(r_0, r_1, r_{\infty}), \] where $r_0 = p + \frac{1}{6} + \delta (l + \frac{1}{2})_{\tilde{+}}$ is as in (\ref{r0}), $r_{\infty} = p+l$ and \[ r_1 = \begin{cases} p + \frac{1}{2} + (l - \frac{1}{2})_{\tilde{+}} &\text{ if } l \geq \frac{1}{2}, \\ p + \frac{1}{2} + \delta (l - \frac{1}{2}) &\text{ if } l < \frac{1}{2}. \end{cases}\] 

For $l > \frac{1}{2}$, one has $r_1 = r_{\infty}$, so $\max(r_0, r_1, r_{\infty}) = \max (p + \frac{1}{6} + \delta(l + \frac{1}{2}), p+l) \geq p+l$, with equality being attained if and only if $\delta \leq (l - \frac{1}{6})/(l + \frac{1}{2})$. Since the right hand side is always strictly larger than $\frac{1}{3}$ in the range of $l$ being considered, choosing $\delta = \frac{1}{3}$ gives $r = p + l$ in this case. Similarly, if $l < -\frac{1}{2}$, then $r = \max(r_0, r_1, r_{\infty}) = p + \frac{1}{6}$ if and only if $\delta \geq \frac{1}{3}/(\frac{1}{2}-l)$. The right hand side being always smaller than $1/3$ in the given range of $l$, choosing any $\frac{1}{3} \leq \delta < \frac{1}{2}$ suffices. For $-\frac12 < l \leq \frac12$, we get
\begin{eqnarray*}
r &= & \min_{\delta \geq \frac{1}{3}} \max \left(p+\delta l +\frac{1-\delta}2,p+\frac16 +\delta(l+\frac12), p+l \right)\\
&= & \min_{\delta \geq \frac{1}{3}} \left[p+\delta l +\max \left(\frac{1-\delta}2,\frac{1+3\delta}6 \right) \right] \\
&= & \min_{\delta \geq \frac{1}{3}} \left[ p+\delta l+\frac{1+3\delta}6 \right] \\ &=& p+\frac{l+1}3.
\end{eqnarray*}
The case $l = -\frac{1}{2}$ is handled similarly. The details are left to the reader. $\qquad \square$ 
     
\subsection{Remark}\label{subsec-remark} Finally, we point out how, for  the geometry $(\D,C_1)$ arising for the normal operator in single source  seismic imaging in the presence of fold caustics,  and more generally for FIOs with folding canonical relations, Thm. \ref{thm-main} follows from the proof for $(\D,C_0)$. From \cite{no}, one knows that the linearized forward scattering operator $F$ belongs to $I^1(C_{sing})$, where, as described in \S\S\ref{subsec-single},  $C_{sing}\subset T^*X\times T^*Y$ is a folding canonical relation, and $F^*F\in I^{2,0}(\Delta, C_1)$. However,  any folding canonical relation  $C_{fold}$ may, by application of canonical transformations $\chi_L,\chi_R$ on the left and right, resp.,   be microlocally put in Melrose-Taylor normal form \cite{meta}, so that we may assume that $C_{fold} \subset C_0\times \D_{T^*\R^{n-2}}$. In composing $C_{fold}\circ C_{fold}$, the $\chi_L,\chi_L^{-1}$ cancel, and a calculation  shows that 
$$\chi_R^{-1}(C_{fold}^t\circ C_{fold})\chi_R\subset \chi_R^{-1}(\D\cup C_1)\chi_R\subset \D_{T^*\R^n}\cup \Bigl[C_0^{\frac14}\times \D_{T^*\R^{n-2}}\Bigr],$$
with $C_0^{\frac{1}{4}} =N^*\{x_2-y_2=\frac14 (x_1-y_1)^3\}'$; on the operator level, for $A\in I^m(C_{fold}), \linebreak B\in I^{m'}(C_{fold})$, one has $B^*A\in I^{m+m',0}(\D,C_1)$ \cite{no,fel}. Since $\chi_R$ and its inverse correspond to unitary FIOs, 
$L^2$-Sobolev estimates for $I^{p,l}(\D,C_1)$ follow from those for $I^{p,l}(\D,C_0^{\frac14}\times \D_{T^*\R^{n-2}})$.
The coefficient $\frac14$ being irrelevant to the analysis, the 
argument above for the $(\D,C_0)$ geometry can be repeated.

\section{Microlocal invertibilty : a counterexample}\label{hilbert-example}
In this section we describe an example related to the question of microlocal invertibility. Given an operator $A$ of the type considered in this paper, e.g., \linebreak$A \in I^{p,l}(\Delta, C_0)$ as in (\ref{Skernel0}), 
that is also elliptic in the  sense that 
\begin{equation} \inf_{x,y,\theta} \frac{|a(x,y;\theta)|}{\langle \theta_2,\theta_1 \rangle^{p+\frac12}
\langle \theta_1 \rangle^{l-\frac12}} \geq c > 0, \label{elliptic} \end{equation} 
we would like to determine whether it is possible to left-invert $A$ microlocally in some appropriate $I^{p',l'}(\Delta, C_0)$ class. More precisely, and in view of Thm. \ref{thm-main}, we pose the following:

{\bf{Question:}} {\em{Let r = r(p,l) denote the regularity exponent obtained in Thm \ref{thm-main}. Given $p,l$, does there exist $p',l' \in \mathbb R$ satisfying \[r(p,l) + r(p',l') = 0 \] and $\epsilon > 0$ depending only on $p,p',l,l'$ such that for every $A \in I^{p,l}(\Delta, C_0)$ that it elliptic in the sense of (\ref{elliptic}), one can find $B \in I^{p',l'}(\Delta, C_0)$ such that 
\begin{equation} E = BA - I \text{ maps $H^s_{\text{comp}}$ boundedly to $H^{s+\epsilon}_{\text{loc}}$ for all $s \in \mathbb R$? } \label{def-smoothing} \end{equation} }}
We will call an operator $E$ satisfying the Sobolev mapping property in (\ref{def-smoothing}) {\em{smoothing of order at least}} $\epsilon$. The following result shows that the answer to this question is no in general. 

\begin{proposition}
Let $\mathcal H$ be the Hilbert transform along the cubic in $\mathbb R^2$ defined in (\ref{eqn-hil}), and let $\psi \in C_0^{\infty}(\mathbb R^2)$. Then the operator $\mathcal H_0 = \psi \mathcal H \psi$ lies in $I^{-\frac{1}{2}, \frac{1}{2}}(\Delta, C_0)$ but there does not exist any $B$ that maps $H^s_{\text{comp}}$ boundedly into $H^s_{\text{loc}}$ for every $s \in \mathbb R$ and for which $B \mathcal H_0 - I$ is smoothing of any positive order.   \label{propn-example}
\end{proposition} 

\subsection{Facts about $\mathcal H$} \label{factsH} Since the operator $\mathcal H$ is translation-invariant, its Fourier transform is a multiplier operator, i.e., there exists a function $m$ on $\mathbb R^2$ such that 
\begin{equation} \widehat{\mathcal H f}(\xi_1, \xi_2) = m(\xi_1, \xi_2) \widehat{f}(\xi_1, \xi_2), \quad \text{ with } \quad m(\xi_1, \xi_2) = \int e^{i(\xi_1 t + \xi_2 t^3)} \frac{dt}{t}. \label{multiplier} \end{equation} 
Properties of  multipliers such as $m$ have been extensively studied in the literature;
see \cite{SW1},\cite{pa} and references there. We note below without proof a few well-known facts about $m$ that are crucial to the proof of the proposition.
\begin{enumerate}[(i)]
\item $m \in L^{\infty}(\mathbb R^2) \cap C^{\infty}(\mathbb R^2 \setminus \{0\})$. \label{hil-item1}
\item $m(\rho \xi_1, \rho^3 \xi_2) = m(\xi_1, \xi_2)$ for all $\rho > 0$ and all $\xi = (\xi_1, \xi_2) \in \mathbb R^2$. In fact, $m(\mu, 1)$ is an antiderivative of the classical Airy function composed with an affine transformation. See   \cite[Ch. 10.4]{AbSt} for details. \label{hil-item2}
\item Figure \ref{mu-plot} contains the graph of the function $\mu \mapsto m(\mu, 1)$. We observe that the graph has a zero, i.e., there exists $\alpha \ne 0$ such that $m(\xi_1, \alpha \xi_1^3) = 0$. In other words, there is a unique cubic curve on which $m$ vanishes.  \label{hil-item3}
\end{enumerate} 
Properties (\ref{hil-item1})-(\ref{hil-item3}) imply that there exists a constant $C > 0$  
\begin{equation} \bigl| m(\xi) \bigr| \leq C \min \left[1, \left|\alpha - \frac{\xi_2}{\xi_1^3}\right| \right]. \label{multiplier-estimate} \end{equation}
\begin{figure}
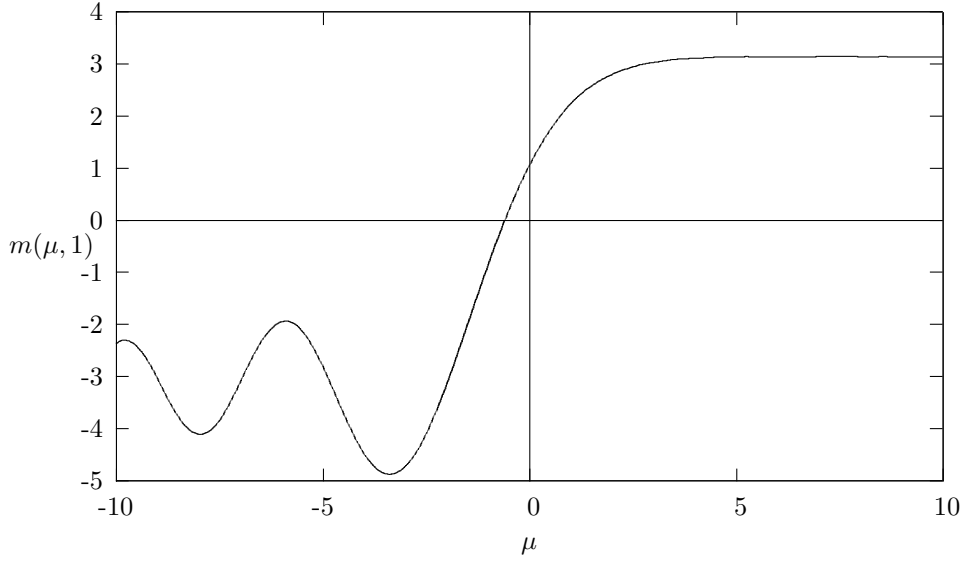

\include{singraph} \caption{Plot of $\mu \mapsto m(\mu,1)$} \label{mu-plot}
\end{figure}
\subsection{Construction of the function $f_0$} The key element of the proof of Proposition \ref{propn-example} is the construction of a function $f_0$ of fixed compact support which lies in $L^2$ and not in $H^s$ for any $s > 0$, but for which $\mathcal Hf \in H^{\epsilon}_{\text{loc}}$ for some $\epsilon > 0$. To this end, let us fix $\chi \in C_0^{\infty}[-1,1]$ such that $\widehat{\chi} \geq 0$ on $\mathbb R$, $\widehat{\chi}(0) > 0$, set 
\[ c_k = k^{-\frac{1}{2}} (1 + \log k)^{-1}, \qquad k \geq 1, \]
and define 
\begin{equation} \widehat{f}_0(\xi) = \sum_{k \geq 1} c_k \widehat{\chi}(\xi_1 - n_k) \widehat{\chi}(\xi_2 - \alpha n_k^3), \label{def-f0} \end{equation}
where $\{n_k \}$ is a fast growing sequence to be determined in the sequel. In particular, we will see that choosing $n_k = 2^k$ works. We collect several facts about $f_0$ in the following sequence of lemmas. 
\begin{lemma}
The function $f_0$ given by (\ref{def-f0}) has compact support in $[-1,1] \times [-1,1]$ and lies in $L^2(\mathbb R^2)$. Aposteriori, this means that the infinite sum in (\ref{def-f0}) converges for almost every $\xi$.  \label{f0-L2} 
\end{lemma} 
\begin{proof}
Since \[f_0(x) = \chi(x_1) \chi(x_2) \sum_{k \geq 1} c_k e^{i(n_k x_1 + \alpha n_k^3 x_2)}, \] the support property of $f_0$ follows easily from that of $\chi$. We proceed to estimate the $L^2$ norm of $f_0$:
\begin{align*}
||f_0||_2^2 = ||\widehat{f}_0||_2^2 &= \sum_{k,k' \geq 1} c_k c_{k'} \int \widehat{\chi}(\xi_1 - n_k) \widehat{\chi}(\xi_2 - \alpha n_k^3) \widehat{\chi}(\xi_1 - n_{k'}) \widehat{\chi}(\xi_2 - \alpha n_{k'}^3) \, d\xi_1 \, d\xi_2 \\ &= \left[ \sum_{k=k'} + \sum_{k \ne k'} \right] =: \mathbf S_1 + \mathbf S_2, \quad  \text{where} \end{align*} 
\begin{align*}  
\mathbf S_1 &= \sum_{k \geq 1} c_{k}^2  \int \left[ \widehat{\chi}(\xi_1 - n_k) \widehat{\chi}(\xi_2 - \alpha n_k^3) \right]^2 \, d\xi \\ &\leq C \sum_{k \geq 1} k^{-1} (1 + \log k)^{-2} < \infty, \text{ and }  \\
\mathbf S_2 &= \sum_{k \ne k'} c_k c_{k'} \int \widehat{\chi}(\xi_1 - n_k) \widehat{\chi}(\xi_2 - \alpha n_k^3) \widehat{\chi}(\xi_1 - n_{k'}) \widehat{\chi}(\xi_2 - \alpha n_{k'}^3) \, d\xi \\ &\leq C \sum_{k \ne k'} c_k c_{k'} \int \widehat{\chi}(\xi_1 - n_k) \widehat{\chi}(\xi_1 - n_{k'}) \, d\xi_1 \\ &\leq C \sum_{k \ne k'} c_k c_{k'} \left[ \int_{|\xi_1 - n_k| \leq \frac{1}{2}|n_k - n_{k'}|} + \int_{|\xi_1 - n_k| \geq \frac{1}{2}|n_k - n_{k'}|} \right] \\ &\leq  C \sum_{k \ne k'} c_k c_{k'} \left[ \int_{|\xi_1 - n_{k'}| \geq \frac{1}{2}|n_k - n_{k'}|} + \int_{|\xi_1 - n_k| \geq \frac{1}{2}|n_k - n_{k'}|} \right] \\ &\leq C \sum_{k \ne k'} c_k c_{k'} \int_{|\xi_1 - n_k| \geq \frac{1}{2}|n_k - n_{k'}|} \widehat{\chi}(\xi_1 - n_k) \, d\xi_1 \\ &\leq C \sum_{k \ne k'} c_k c_{k'} \int \left[ 1 + |\xi_1 - n_k| \right]^{-4} \, d\xi_1 \\ &\leq C \sum_{k \ne k'} c_k c_{k'} (1 + |n_k - n_{k'}|)^{-2} \\ &\leq C \sum_{k \ne k'} \frac{c_k^2}{(1 + |n_k - n_{k'}|)^2} \leq C \sum_k c_k^2 < \infty. 
\end{align*}  
In the computation above, the second step in the estimation of $\mathbf S_2$ follows from the bound 
\[ \sup_{k,k'} \int \widehat{\chi}(\xi_2 - \alpha n_k^3) \widehat{\chi}(\xi_2 - \alpha n_{k'}^3) \, d\xi_2 \leq C < \infty, \] 
while the fourth step uses the inequality $|\xi_1 - n_{k'}| \geq |n_k - n_{k'}| - |\xi_1 - n_k| \geq \frac{1}{2}|n_k - n_{k'}|$ if $|\xi_1 - n_k| \leq \frac{1}{2}|n_k - n_{k'}|$. This completes the proof of Lemma \ref{f0-L2}. 
\end{proof}   
\begin{lemma}
The function $f_0$ given by (\ref{def-f0}) does not lie in $H^s_{\text{loc}}$ for any $s > 0$. \label{f0-notinHs}
\end{lemma} 
\begin{proof}
Let $\psi$ be any smooth function of compact support. Without loss of generality, we may assume that $\widehat{\psi} \geq 0$, with $\widehat{\psi}(0) > 0$.  We have to show that $\psi f_0 \not \in H^s$ for any $s > 0$.  

Let us now choose small constants $c_0$, $\epsilon_0 > 0$ such that \begin{equation} \inf \{\widehat{\psi}(\xi) : |\xi| \leq \epsilon_0 \} \geq c_0, \quad \inf\{\widehat{\chi}(\xi_1) : |\xi_1| \leq \epsilon_0 \} \geq c_0. \label{psi-chi-lowerbound} \end{equation}
We compute 
 \begin{align*}
||\psi f_0||_{H^s}^2 &= \int (1 + |\xi|^2)^{s} \bigl| \widehat{\psi} \ast \widehat{f}_0(\xi) \bigr|^2 \, d\xi \\ &= \int (1 + |\xi|^2)^{s} \left[\int \widehat{\psi}(\xi - \eta) \widehat{f}_0(\eta)\, d\eta \right] \left[\int \widehat{\psi}(\xi - \eta') \widehat{f}_0(\eta') \, d\eta' \right] \, d\xi \\ &\geq c_0^2 \iiint_{\mathcal D} (1 + |\xi|^2)^{s} \widehat{f}_0(\eta) \widehat{f}_0(\eta') \, d\eta \, d\eta' \, d\xi, 
\end{align*} 
where at the last step we have replaced the domain of $(\xi, \eta, \eta')$ integration by the subset 
\[ \mathcal D := \left\{(\xi, \eta, \eta') : |\xi - \eta| \leq \epsilon_0, |\xi - \eta'| \leq \epsilon_0 \right\}, \] and then used the first inequality in (\ref{psi-chi-lowerbound}). Substituting the expression for $\widehat{f}_0$ from (\ref{def-f0}) in the last integrand and replacing the infinite sum (of non-negative summands) in $k, k'$ by the smaller subsum over $k = k'$ we obtain 
\begin{align*}
||\psi f_0||_{H^s}^2 &\geq c_0 \sum_{k,k' \geq 1} c_{k} c_{k'} \iiint_{\mathcal D} (1 + |\xi|^2)^{s} \widehat{\chi}(\eta_1 - n_k) \widehat{\chi}(\eta_2 - \alpha n_k^3) \\ &\hskip2in \times \widehat{\chi}(\eta_1' - n_{k'}) \widehat{\chi}(\eta_2' - \alpha n_{k'}^3) \, d\eta \, d\eta' \, d\xi \\ &\geq c_0^2 \sum_{k \geq 1} c_k^2 \iiint_{\mathcal D} (1 + |\xi|^2)^{s} \widehat{\chi}(\eta_1 - n_k) \widehat{\chi}(\eta_2 - \alpha n_k^3) \\ &\hskip2in \times \widehat{\chi}(\eta_1' - n_k) \widehat{\chi}(\eta_2' - \alpha n_k^3) \, d\eta \, d\eta' \, d\xi \\ &\geq c_0^6 \sum_{k \geq 1} c_k^2 \iiint_{(\xi, \eta, \eta') \in \mathcal D_k \times \mathcal D_k \times \mathcal D_k} \left(\frac{1}{2}(1 + |\xi|)^2 \right)^{s} \, d\eta \, d\eta' \, d\xi,  
\end{align*} 
where \[\mathcal D_k := \left\{\xi : |\xi - (n_k, \alpha n_k^3)| < \frac{\epsilon_0}{2} \right\}, \quad \text{so that } \quad \mathcal D_k \times \mathcal D_k \times \mathcal D_k \subseteq \mathcal D, \]
and the lower bound for the integrand in the last step follows from the second inequality in (\ref{psi-chi-lowerbound}). On $\mathcal D_k$, the variable $\xi$ satisfies the lower bound $1 + |\xi| \geq 1 + |\alpha n_k^3| - |\xi_2 - \alpha n_k^3| \geq 1 - \epsilon_0 + |\alpha n_k^3| \geq |\alpha n_k^3|$, so that there exists a constant $c > 0$ for which  
\[ ||\psi f_0||_{H^s}^2 \geq c \sum_{k \geq 1} c_k^2 n_k^{6s}. \]
The sum on the right hand side above diverges to $\infty$ for any $s > 0$ if $n_k$ increases sufficiently fast, for instance exponentially. This completes the proof of the Lemma \ref{f0-notinHs}.   
\end{proof} 
\begin{lemma}
There exists $s_0 > 0$ such that $\mathcal Hf_0 \in H^{s_0}$. \label{Hf0-inHs}  
\end{lemma} 
\begin{proof}
It follows from (\ref{multiplier}) and (\ref{def-f0}) that 
\begin{align*}
||\mathcal Hf_0||_{H^s}^2 &= \int \left| \widehat{\mathcal H f_0}(\xi) \right|^2 (1 + |\xi|^2)^{s} \, d\xi \\ &= \int |m(\xi)|^2  |\widehat{f}_0(\xi)|^2 (1 + |\xi|^2)^{s} \, d\xi \\ &= \sum_{k,k'} c_k c_{k'} \int \bigl| m(\xi) \bigr|^2 (1 + |\xi|^2)^s \widehat{\chi}(\xi_1 - n_k) \widehat{\chi}(\xi_2 - \alpha n_k^3) \\ &\hskip1.5in \times \widehat{\chi}(\xi_1 - n_{k'}) \widehat{\chi}(\xi_2 - \alpha n_{k'}^3) \, d\xi \\ 
&= \left[ \sum_{k = k'} + \sum_{k \ne k'} \right] =: \mathbf T_1 + \mathbf T_2. 
\end{align*}  
We now proceed to estimate the two sums separately. For the first, we write 
\begin{equation}
\begin{aligned} 
\mathbf T_1 &= \sum_{k \geq 1} c_k^2 \int \bigl|m(\xi) \bigr|^2 \widehat{\chi}(\xi_1 - n_k)^2 \widehat{\chi}(\xi_2 - \alpha n_k^3)^2 (1 + |\xi|^2)^s \, d\xi \\ &= \sum_{k \geq 1} c_k^2 \left[ \int_{\mathcal E_k}  + \int_{\mathcal E_k^c}\right] =: \mathbf T_{11} + \mathbf T_{12}, \quad \text{ where } \\
\mathcal E_k &= \left\{\xi : |\xi_1 - n_k| \leq n_k^{\kappa}, \quad |\xi_2 - \alpha n_k^3| \leq n_k^{\kappa} \right\}.  
\end{aligned} \label{T1} \end{equation} 
Here $\kappa$ is a fixed small positive constant (in fact, any $\kappa < \frac{1}{3}$ will work). The multiplier estimate (\ref{multiplier-estimate}) yields 
\begin{equation} \begin{aligned}
\mathbf T_{11} &\leq C \sum_{k \geq 1} c_k^2 \int_{\mathcal E_k} \left|\alpha - \frac{\xi_2}{\xi_1^3}  \right|^2 (1 + |\xi|)^{2s} \, d\xi \\ &\leq C \sum_{k \geq 1} c_k^2 \sup_{\xi \in \mathcal E_k} \left[ \left(1 + |\xi| \right)^{2s} \frac{|\xi_2 - \alpha \xi_1^3|^2}{|\xi_1|^6}\right] |\mathcal E_k|.  
\end{aligned} \label{T11} \end{equation} 
But on $\mathcal E_k$,
\begin{equation} 
\begin{aligned}
\bigl|\xi_2 - \alpha \xi_1^3 \bigr| &= \bigl|\xi_2 - \alpha \bigl((\xi_1 - n_k) + n_k \bigr)^3\bigr| \\ &= \bigl|\xi_2 - \alpha n_k^3 \bigr| + |\alpha| |\xi_1 - n_k|^3 \\ &\hskip1.5in +3|\alpha| n_k |\xi_1 - n_k|^2 + 3 |\alpha| n_k^2 |\xi_1 - n_k| \\ &\leq n_k^{\kappa} + |\alpha| n_k^{3 \kappa} + 3|\alpha| n_k^{1 + 2 \kappa} + 3|\alpha| n_k^{2 + \kappa} \\ &\leq C n_k^{2 + \kappa}, \end{aligned} \label{estimatea} \end{equation} 
while \begin{equation}  
|\xi_1| \geq n_k - |\xi_1 - n_k| \geq n_k - n_k^{\kappa} \geq \frac{1}{2} n_k, \label{estimateb} \end{equation} and 
\begin{equation} \begin{aligned} 1 + |\xi| &\leq 1 + |\xi_1| + |\xi_2| \\ &\leq 1 + n_k + |\alpha n_k^3| + |\xi_1 - n_k| + |\xi_2 - \alpha n_k^3| \\ &\leq 1 + n_k + |\alpha n_k^3| + 2 n_k^{\kappa} \leq C n_k^3. \end{aligned} \label{estimatec}  \end{equation} 
Combining estimates (\ref{estimatea}), (\ref{estimateb}) and (\ref{estimatec}), and using the fact that $|\mathcal E_k| = n_k^{2 \kappa}$, we arrive at the following bound for $\mathbf T_{11}$:
\[ \mathbf T_{11} \leq C \sum_{k \geq 1} c_k^2 n_k^{6s} \frac{(n_k^{2+\kappa})^2}{n_k^6} n_k^{2 \kappa} \leq C \sum_{k \geq 1} c_k^2 n_k^{4 \kappa + 6s - 2} < \infty \]   
if $s > 0$ is chosen small enough so that $4 \kappa + 6s < 2$. 

We will now estimate $\mathbf T_{12}$. For any $N \geq 1$, 
\begin{align*}
\mathbf T_{12} &= \sum_{k \geq 1} c_k^2 \int_{\mathcal E_k^c} \bigl|m(\xi) \bigr|^2 \widehat{\chi}(\xi_1 - n_k)^2 \widehat{\chi}(\xi_2 - \alpha n_k^3)^2 \bigl(1 + |\xi|^2 \bigr)^{s} \, d\xi \\ &\leq C_N \sum_{k \geq 1} c_k^2 \int_{\mathcal E_k^c} \frac{(1 + |\xi_1|)^{2s}(1 + |\xi_2|)^{2s} \, d\xi}{(1 + |\xi_1 - n_k|)^{2N} (1 + |\xi_2 - \alpha n_k^{3}|)^{2N}} \\ &\leq C_N \sum_{k \geq 1} c_k^2 \int_{\mathcal E_k^c} \frac{(1 + n_k + |\xi_1 - n_k|)^{2s}(1 + |\alpha n_k^3| + |\xi_2 - \alpha n_k^3|)^{2s} \, d\xi}{(1 + |\xi_1 - n_k|)^{2N} (1 + |\xi_2 - \alpha n_k^{3}|)^{2N}} \\ &\leq C_N \sum_{k \geq 1} c_k^2 (1 + n_k)^{2s} (1 + |\alpha n_k^3|)^{2s} \int_{\mathcal E_k^c} \frac{d \xi}{\left[(1 + |\xi_1 - n_k|)(1 + |\xi_2 - \alpha n_k^3|)\right]^{2N - 2s}} \\ &\leq C_N \sum_{k \geq 1} c_k^2 n_k^{8s - \kappa N} \int \bigl( 1 + |\xi_1 - n_k| \bigr)^{-2} \bigl(1 + |\xi_2 - \alpha n_k^3| \bigr)^{-2} \, d\xi < \infty.    
\end{align*} 
The last step follows from the fact that on $\mathcal E_k^c$ at least one of the quantities $|\xi_1 - n_k|$ or $|\xi_2 - \alpha n_k^3|$ must be $\geq n_k^{\kappa}$. For the sum in the last step to converge for a given choice of $s$ and $\kappa$, we must choose $N$ large enough so that $2s - N < -2$ and $8s - \kappa N < 0$. 

Finally, we turn to the estimation of $\mathbf T_2$. This will be done almost exactly in the same way as $\mathbf T_1$ but this time keeping in mind the almost orthogonality of the summands. More precisely, 
\begin{align*}
\mathbf T_2 &= \sum_{k \ne k'} c_k c_{k'} \int \bigl|m(\xi)\bigr|^2 (1 + |\xi|^2)^s  \widehat{\chi}(\xi_1 - n_k) \widehat{\chi}(\xi_2 - \alpha n_k^3) \\ &\hskip2.8in \times \widehat{\chi}(\xi_1 - n_{k'}) \widehat{\chi}(\xi_2 - \alpha n_{k'}^3) \,  d\xi \\ &= \sum_{k \ne k'} c_k c_{k'} \left[ \int_{|\xi_1 - n_k| \leq \frac{1}{2}|n_k - n_{k'}|} + \int_{|\xi_1 - n_k| \geq \frac{1}{2}|n_k - n_{k'}|}\right] \\ & \leq  \sum_{k \ne k'} c_k c_{k'} \left[ \int_{|\xi_1 - n_{k'}| \geq \frac{1}{2}|n_k - n_{k'}|} + \int_{|\xi_1 - n_k| \geq \frac{1}{2}|n_k - n_{k'}|} \right] \\ &\leq C \sum_{k \ne k'} c_{k} c_{k'} \int_{|\xi_1 - n_{k'}| \geq \frac{1}{2}|n_k - n_{k'}|} \bigl|m(\xi)\bigr|^2 (1 + |\xi|^2)^s  \widehat{\chi}(\xi_1 - n_k) \widehat{\chi}(\xi_2 - \alpha n_k^3) \\ &\hskip2.8in \times \widehat{\chi}(\xi_1 - n_{k'}) \widehat{\chi}(\xi_2 - \alpha n_{k'}^3) \,  d\xi
\\ & \leq C \sum_{k \ne k'} \frac{c_k c_{k'}}{(1 + |n_k - n_{k'}|)^2} \int \bigl| m(\xi) \bigr|^2 (1 + |\xi|^2)^{s} \widehat{\chi}(\xi_1 - n_k) \widehat{\chi}(\xi_2 - \alpha n_k^3) \,  d\xi.   
\end{align*} 
We observe that the integrand in the last step above is of the same form as the one occurring in the expression for $\mathbf T_1$ (see the first line in (\ref{T1})); the only distinction is that the function $\widehat{\chi}^2$ in $\mathbf T_1$ has been replaced by $\widehat{\chi}$, which makes no difference to the estimation process. We leave the reader to verify that the same set of arguments that was used to estimate $\mathbf T_1$ now yields a number $r = \min (2 - 4 \kappa - 6s, \kappa N - 8s) > 0$ such that 
\[ \int \bigl| m(\xi)\bigr|^2 \widehat{\chi}(\xi_1 - n_k) \widehat{\chi}(\xi_2 - \alpha n_{k}^3) (1 + |\xi|)^{2s} \, d\xi \leq C n_k^{-r},   \] 
so that 
\begin{align*}
\mathbf T_2 &\leq C \sum_{k \ne k'} \frac{c_k c_{k'} n_k^{-r}}{(1 + |n_k - n_{k'}|)^2} \\ & \leq C \left[\sum_{k \ne k'} \frac{c_{k}^2 n_{k}^{-2r}}{1 + |n_k - n_{k'}|)^2} \right]^{\frac{1}{2}} \left[\sum_{k \ne k'} \frac{c_{k'}^2}{1 + |n_k - n_{k'}|)^2} \right]^{\frac{1}{2}} < \infty. 
\end{align*} 
This completes the proof of Lemma \ref{Hf0-inHs}. 
\end{proof} 
\subsection{Proof of Proposition \ref{propn-example}} It is easy to check that $\mathcal H_0 \in I^{-\frac{1}{2}, \frac{1}{2}}(\Delta, C_0)$, and we leave this to the reader. Without loss of generality, we may assume that the function $\psi$ in the statement of the proposition is identically 1 on $[-1,1]^2$ and vanishes outside $[-2,2]^2$. Let us suppose that there exists an operator $B$ that maps $H^{s}_{\text{comp}}$ boundedly to $H^s_{\text{loc}}$ for all $s \in \mathbb R$ and such that $B \mathcal H_0 = I + E$, where $E$ is a smoothing operator of some order $\geq s_0 > 0$. Then $f_0 = B \mathcal H_0 f_0 - Ef_0$, where $f_0$ is the function in (\ref{def-f0}). Since supp$(f_0) \subseteq [-1,1]^2$ by Lemma \ref{f0-L2}, it follows from Lemma \ref{Hf0-inHs} that $\mathcal H(\psi f_0) = \mathcal Hf_0 \in H^{s_1}$ for some small $s_1 > 0$, so $B\mathcal H_0 f_0 \in H^{s_1}_{\text{loc}}$. On the other hand, $Ef_0 \in H^{s_0}_{\text{loc}}$ since $f_0 \in L^2$ by Lemma \ref{f0-L2}. This implies that $f_0 = B \mathcal H_0 f_0 - Ef_0 \in H^{s_2}_{\text{loc}}$ where $s_2 = \min(s_0, s_1) > 0$, thereby contradicting Lemma \ref{f0-notinHs}.  
\hfill $\square$

\vfil\eject

\vskip.2in

\noindent{\sc School of of Mathematical Sciences}

\noindent{\sc Rochester  Institute of Technology}

\noindent{\sc  Rochester, NY 14623}

\noindent{\tt{rxfsma@rit.edu}}

\vskip.2in

\noindent{\sc Department of Mathematics}

\noindent{\sc University of Rochester}

\noindent{\sc Rochester, NY 14627}

\noindent{\tt{allan@math.rochester.edu}}

\vskip.2in

\noindent{\sc Department of Mathematics}

\noindent{\sc University of British Columbia}

\noindent{\sc Vancouver, B.C., Canada V6T 1Z2}

\noindent{\tt{malabika@math.ubc.ca}}

\end{document}